\documentclass[a4paper]{article}
\usepackage{amsmath,amsthm,amsfonts,thmtools}
\usepackage{amssymb}
\usepackage{graphicx}
\usepackage{subfigure}
\usepackage{euscript}
\usepackage{units}
\usepackage{enumerate}
\usepackage{enumitem}
\usepackage{xcolor}
\usepackage{bbm}

\numberwithin{equation}{section}
\newtheorem{theorem}{Theorem}[section]
\newtheorem{maintheorem}[theorem]{Main Theorem}
\newtheorem{lemma}[theorem]{Lemma}

\newtheorem{proposition}[theorem]{Proposition}

\theoremstyle{definition}

\newtheoremstyle{myremarkstyle}{}{}{}{}{\bfseries}{.}{ }{}
\theoremstyle{myremarkstyle}
\declaretheorem[name=Remark,qed=$\blacksquare$,numberlike=theorem]{remark}
\declaretheorem[name=Remarks,qed=$\blacksquare$,numberlike=theorem]{remarks}

\declaretheorem[name=Definition,qed=$\blacksquare$,numberlike=theorem]{definition}
\declaretheorem[name=Notation,qed=$\blacksquare$,numberlike=theorem]{notation}

\makeatletter
\newcommand*{\intavg}{%
  \mint@l{-}{}%
}
\newcommand*{\mint@l}[2]{%
  \@ifnextchar\limits{%
    \mint@l{#1}%
  }{%
    \@ifnextchar\nolimits{%
      \mint@l{#1}%
    }{%
      \@ifnextchar\displaylimits{%
        \mint@l{#1}%
      }{%
        \mint@s{#2}{#1}%
      }%
    }%
  }%
}
\newcommand*{\mint@s}[2]{%
  \@ifnextchar_{%
    \mint@sub{#1}{#2}%
  }{%
    \@ifnextchar^{%
      \mint@sup{#1}{#2}%
    }{%
      \mint@{#1}{#2}{}{}%
    }%
  }%
}
\def\mint@sub#1#2_#3{%
  \@ifnextchar^{%
    \mint@sub@sup{#1}{#2}{#3}%
  }{%
    \mint@{#1}{#2}{#3}{}%
  }%
}
\def\mint@sup#1#2^#3{%
  \@ifnextchar_{%
    \mint@sub@sup{#1}{#2}{#3}%
  }{%
    \mint@{#1}{#2}{}{#3}%
  }%
}
\def\mint@sub@sup#1#2#3^#4{%
  \mint@{#1}{#2}{#3}{#4}%
}
\def\mint@sup@sub#1#2#3_#4{%
  \mint@{#1}{#2}{#4}{#3}%
}
\newcommand*{\mint@}[4]{%
  \mathop{}%
  \mkern-\thinmuskip
  \mathchoice{%
    \mint@@{#1}{#2}{#3}{#4}%
        \displaystyle\textstyle\scriptstyle
  }{%
    \mint@@{#1}{#2}{#3}{#4}%
        \textstyle\scriptstyle\scriptstyle
  }{%
    \mint@@{#1}{#2}{#3}{#4}%
        \scriptstyle\scriptscriptstyle\scriptscriptstyle
  }{%
    \mint@@{#1}{#2}{#3}{#4}%
        \scriptscriptstyle\scriptscriptstyle\scriptscriptstyle
  }%
  \mkern-\thinmuskip
  \int#1%
  \ifx\\#3\\\else_{#3}\fi
  \ifx\\#4\\\else^{#4}\fi  
}
\newcommand*{\mint@@}[7]{%
  \begingroup
    \sbox0{$#5\int\m@th$}%
    \sbox2{$#5\int_{}\m@th$}%
    \dimen2=\wd0 %
    \let\mint@limits=#1\relax
    \ifx\mint@limits\relax
      \sbox4{$#5\int_{\kern1sp}^{\kern1sp}\m@th$}%
      \ifdim\wd4>\wd2 %
        \let\mint@limits=\nolimits
      \else
        \let\mint@limits=\limits
      \fi
    \fi
    \ifx\mint@limits\displaylimits
      \ifx#5\displaystyle
        \let\mint@limits=\limits
      \fi
    \fi
    \ifx\mint@limits\limits
      \sbox0{$#7#3\m@th$}%
      \sbox2{$#7#4\m@th$}%
      \ifdim\wd0>\dimen2 %
        \dimen2=\wd0 %
      \fi
      \ifdim\wd2>\dimen2 %
        \dimen2=\wd2 %
      \fi
    \fi
    \rlap{%
      $#5%
        \vcenter{%
          \hbox to\dimen2{%
            \hss
            $#6{#2}\m@th$%
            \hss
          }%
        }%
      $%
    }%
  \endgroup
}

\DeclareMathOperator*{\esssup}{ess\ sup}
\DeclareMathOperator{\supp}{supp}

\DeclareMathOperator*{\wlim}{w-lim}
\DeclareMathOperator*{\wslim}{w*-lim}
\DeclareMathOperator{\Lip}{Lip}

\newcommand{\wto}{\rightharpoonup}
\newcommand{\wsto}{\overset{*}{\wto}}

\renewcommand{\geq}{\geqslant}
\renewcommand{\leq}{\leqslant}

\newcommand{\eps} {\varepsilon}

\renewcommand{\epsilon}{\varepsilon}
\renewcommand{\phi}{\varphi}

\newcommand{\Lone}{\EuScript{H}}

\newcommand{\A}{\mathcal{A}}
\newcommand{\R}{\mathbb{R}}
\newcommand{\N}{\mathbb{N}}

\newcommand{\F}{\EuScript{F}}
\newcommand{\M}{\EuScript{M}}

\newcommand{\U}{U}		% Phase space

\newcommand{\Borel}{\EuScript{B}}

\newcommand{\Prob}{\EuScript{P}}
\newcommand{\Lp}{\EuScript{L}}

\newcommand{\ind}{\mathbbm{1}} 
\newcommand{\Cyl}{\mathit{Cyl}}

\newcommand{\define}[1]{\textbf{#1}}

\newcommand{\ip}[2]{\bigl\langle #1,\, #2\bigr\rangle}	% Pairing of measure and function
\newcommand{\lip}[2]{#1(#2)}				% Pairing of L2 functions

\newcommand{\sgn}{{\rm sgn}}

\begin{document}

\title{Statistical solutions of hyperbolic conservation laws I: Foundations}
\author{U. S. Fjordholm\thanks{Department of Mathematical Sciences, Norwegian University of Science and Technology, Trondheim, N-7491, Norway.}, \, S. Lanthaler\thanks{Swiss plasma center, SB SPC-TH, PPB 313 (B\^atiment PPB), Station 13 CH-1015 Lausanne, Switzerland} ,\, and S. Mishra \thanks{Seminar for Applied Mathematics, ETH Z\"urich, R\"amistrasse 101, Z\"urich, Switzerland.}}

\maketitle
\begin{abstract}
We seek to define statistical solutions of hyperbolic systems of conservation laws as time-parametrized probability measures on $p$-integrable functions. To do so, we prove the equivalence between probability measures on $L^p$ spaces and infinite families of \textit{correlation measures}. Each member of this family, termed a \textit{correlation marginal}, is a Young measure on a finite-dimensional tensor product domain and provides information about multi-point correlations of the underlying integrable functions. We also prove that any probability measure on a $L^p$ space is uniquely determined by certain moments (correlation functions) of the equivalent correlation measure.  

We utilize this equivalence to define statistical solutions of multi-dimensional conservation laws in terms of an infinite set of equations, each evolving a moment of the correlation marginal. These evolution equations can be interpreted as augmenting entropy measure-valued solutions, with additional information about the evolution of all possible multi-point correlation functions. Our concept of statistical solutions can accommodate uncertain initial data as well as possibly non-atomic solutions even for atomic initial data. 

For multi-dimensional scalar conservation laws we impose additional entropy conditions and prove that the resulting \textit{entropy statistical solutions} exist, are unique and are stable with respect to the $1$-Wasserstein metric on probability measures on $L^1$.

\end{abstract}

\section{Introduction}
Systems of conservation laws are nonlinear partial differential equations of the generic form
\begin{subequations}\label{eq:clcauchy}
\begin{equation}\label{eq:cl}
\partial_t u + \nabla_x\cdot f(u) = 0
\end{equation}
\begin{equation}\label{eq:clinitial}
u(x,0) = \bar{u}(x).
\end{equation}
\end{subequations}
Here, the unknown $u = u(x,t) : \R^d\times\R_+ \to \R^N$ is the vector of \emph{conserved variables} and $f = (f^1, \dots, f^d) : \R^N \to \R^{N\times d}$ is the \emph{flux function}. We denote $\R_+ := [0,\infty)$. The system is termed \emph{hyperbolic} if the flux Jacobian matrix has real eigenvalues \cite{DAF1}. Here and in the remainder, quantities with a bar (like $\bar{u}$) denote prescribed initial data.

Hyperbolic systems of conservation laws arise in a wide variety of models in physics and engineering. Prototypical examples include the compressible Euler equations of gas dynamics, the shallow water equations of oceanography, the magneto-hydrodynamics (MHD) equations of plasma physics and the equations of nonlinear elasticity \cite{DAF1}.

It is well known that solutions of \eqref{eq:clcauchy} can form discontinuities such as \emph{shock waves}, even for smooth initial data $\bar{u}$. Hence, solutions of systems of conservation laws \eqref{eq:clcauchy} are sought in the sense of distributions. These \emph{weak solutions} are not necessarily unique. They need to be augmented with additional admissibility criteria, often termed \emph{entropy conditions}, to single out the physically relevant solution. \emph{Entropy solutions} are widely regarded as the appropriate solution paradigm for systems of conservation laws \cite{DAF1}. 

Global well-posedness (existence, uniqueness and continuous dependence on initial data) of entropy solutions of scalar conservation laws ($N = 1$ in \eqref{eq:clcauchy}), was established in the pioneering work of Kruzkhov \cite{Kruz1}. For one-dimensional systems ($d=1$, $N>1$ in \eqref{eq:clcauchy}), global existence, under the assumption of small initial total variation, was shown by Glimm in \cite{GL1} and by Bianchini and Bressan in \cite{BB1}. Uniqueness and stability of entropy solutions for one-dimensional systems has also been shown; see \cite{BRE1} and references therein. 

Although existence results have been obtained for some very specific examples of multi-dimensional systems (see \cite{SBG1} and references therein), there are \emph{no global existence results} for any generic class of multi-dimensional systems. In fact, De Lellis, Sz\'ekelyhidi et al.\ have recently been able to construct \emph{infinitely many} entropy solutions for prototypical multi-dimensional systems such as the Euler equations for polytropic gas dynamics (see \cite{DLS1,CDL3} and references therein). Their construction involves a novel iterative procedure where oscillations at smaller and smaller scales are successively added to suitably constructed sub-solutions of \eqref{eq:clcauchy}.

Given the lack of global existence and uniqueness results for entropy solutions of multi-dimensional systems of conservation laws, it is natural to seek alternative solution paradigms. One option, advocated for instance in \cite{BTW1}, is to augment entropy solutions with further admissibility criteria, such as the vanishing viscosity limit, in order to rule out ``unphysical'' solutions. However, given the difficulties of obtaining existence results for the weaker concept of entropy solutions, it is unclear if such a narrowing of the solution concept would lead to any meaningful global existence results. 

The other alternative is to extend the solution concept beyond entropy solutions (integrable functions) and seek possibly even weaker notions of solutions of \eqref{eq:clcauchy}, together with suitable admissibility criteria to constrain these solutions and enforce uniqueness. A recent paper \cite{FKMT15} advocates such an approach. Based on the extensive numerical simulations reported in \cite{FKMT15} (see also \cite{GL2}), the authors observe that approximate solutions of \eqref{eq:clcauchy} can feature oscillations at smaller and smaller scales as mesh is refined. Given this fact, they postulate that \emph{entropy measure-valued solutions} may serve as an appropriate solution paradigm for systems of conservation laws in several space dimensions, particularly in characterizing limits of (numerical) approximations. 

Measure-valued solutions, originally proposed by DiPerna in \cite{DiP85} (see also \cite{DM87}), are space-time-parametrized probability measures, or \emph{Young measures}, defined on the phase space $\R^N$ of \eqref{eq:clcauchy}. In defining entropy measure-valued solutions, one requires consistency of certain functionals of this Young measure with the initial data, with the weak (distributional) form of \eqref{eq:clcauchy}, and with a suitable (dissipative) form of the entropy conditions (see also \cite{DST12}). 

In recent papers \cite{FKMT15,FMTacta} (see also \cite{GGW1}), the authors were able to prove (global in time) existence of entropy measure-valued solutions for a very large class of systems of conservation laws, namely those endowed with a strictly convex entropy function, by showing convergence of numerical approximations of \eqref{eq:clcauchy} based on a Monte Carlo algorithm. Numerical experiments presented in these papers suggest that the measure-valued solution may be \emph{non-atomic}, even when the initial data is atomic, i.e.\ a Dirac Young measure concentrated on an integrable function. The computed measure-valued solutions were observed to be stable with respect to the choice of numerical method and with respect to perturbations of initial data. 

However, one can readily construct counter-examples to uniqueness of these \emph{entropy measure-valued solutions}. In particular, if the initial data is \emph{non-atomic} then infinitely many entropy measure-valued solutions can be constructed, even for scalar conservation laws (see \cite{Sch89,FKMT15}). This lack of uniqueness, even for the scalar case, can be attributed to the fact that only certain functionals of the measure-valued solution (essentially the mean and the second moment) are required to be consistent with the initial data, the evolution equation \eqref{eq:clcauchy} and the entropy conditions. Since the the mean and the second moment uniquely specifies a measure \emph{only} when the measure is atomic, one cannot expect uniqueness for generic (non-atomic) measure-valued solutions as considered in \cite{FKMT15}.

On the other hand, numerical experiments presented in \cite{FKMT15} clearly suggest that one has to deal with non-atomic, ``uncertain'' measure-valued solutions of multi-dimensional systems of conservation laws, even when the initial data is atomic. In a wide variety of applications, even the initial data can be non-atomic, carrying some uncertainty due to e.g.\ measurement errors. These measurements are inherently uncertain and can only be specified probabilistically, and this uncertainty inevitably propagates into the solution. The modeling, analysis and numerical approximation of uncertain solutions, given uncertain inputs (such as the initial data), falls under the rubric of \emph{uncertainty quantification}; see \cite{lncse:uq} and reference therein for an extensive discussion of the very large body of recent research activity on uncertainty quantification for systems of conservation laws. Thus, in general, one has to deal with the possibility that physically relevant measure-valued solutions are non-atomic.

Given these considerations, we seek to find a solution framework that can deal with non-atomic measure-valued solutions of multi-dimensional systems of conservation laws, and can provide further constraints on these measure-valued solutions in order to enforce uniqueness and stability of the resulting solution concept.

A natural choice for such a solution framework is the notion of \emph{statistical solutions} that was first proposed by Foia\c s in \cite{Foi72,Foi73} (see also \cite{FMRT}) in the context of the incompressible Navier--Stokes equations of fluid dynamics. As envisaged by Foia\c s and co-workers, statistical solutions of the Navier--Stokes equations are time-parametrized probability measures on a given infinite-dimensional function space (divergence-free $L^2$ functions in the context of the Navier--Stokes equations). This family of measures has to satisfy either a suitable \emph{infinite-dimensional Liouville equation} that governs the time evolution of a class of functionals in a manner consistent with the Navier--Stokes dynamics, or equivalently, satisfy a \emph{Hopf equation}, where the time-evolution of the characteristic functional of the probability measure (on $L^2$) is prescribed. Both formulations result in evolution equations in infinite-dimensional spaces. A detailed account of statistical solutions in the sense of Foia\c s, and their relation to the description of turbulent incompressible flows, can be found in \cite{FMRT} and references therein.

However, it is far from straightforward to adapt the notion of statistical solutions to the context of systems of conservation laws. There seems to be at least three main difficulties in this regard. First, statistical solutions as defined in \cite{Foi72,Foi73,FMRT} are well suited to problems with viscosity, as they require some regularity of the underlying functions in order to define the infinite-dimensional Liouville or Hopf equations. It is unclear how to extend them to inviscid problems such as systems of conservation laws where solutions are generally discontinuous. Attempts to do so have been made in \cite{CD1,CD2,Bert} (see also \cite{IL1, Pan1}) for the special case of the one-dimensional inviscid Burgers equation. The corresponding statistical solutions are probability measures on the space of \emph{distributions}, and the infinite-dimensional Hopf equation is well-defined by using compactly supported infinitely differentiable test functions. Although existence results for such statistical solutions of the inviscid Burgers equation have been obtained in the class of Levy processes with negative jumps, it is not possible to obtain uniqueness of these statistical solutions, even for the inviscid Burgers equation, in the class of probability measures on spaces as large as the space of distributions. 

The second difficulty with statistical solutions in the sense of Foia\c s, lies in the fact that the Liouville or Hopf equations are evolution equations on infinite-dimensional function spaces. This makes the interpretation and computation of statistical solutions very hard for viscous problems, and the solution concept is not easily amenable to extension to inviscid PDEs such as systems of conservation laws. Furthermore, probability measures on function spaces preclude a local (in space) description of the resulting solution, as it is unclear how to interpret statistical information at specific points (or collection of points) in space.

Finally, given our original motivation in constraining measure-valued solutions to recover uniqueness in the non-atomic case, the relationship between statistical solutions and measure-valued solutions is far from clear. The only known results are presented in \cite{Chae1,Chae2} where a sequence of statistical solutions of the incompressible Navier--Stokes equations is shown to converge to a measure-valued solution of the incompressible Euler equations, as defined in \cite{DM87}, when the viscosity vanishes. However, we are interested in investigating the more abstract question of the relationship between probability measures on function spaces (statistical solutions), and Young measures that represent one-point statistics (measure-valued solutions), with the aim of imposing further constraints on measure-valued solutions to enforce uniqueness.

With this background, the \emph{first aim} of the current paper is to provide a novel representation of a probability measure on an infinite-dimensional function space (to be specific, $L^p$ space) in terms of an \emph{infinite hierarchy} of Young measures called a \emph{correlation measure}, defined on tensor products of the (finite-dimensional) spatial domain. Each member of this hierarchy of measures, termed a \emph{correlation marginal}, represents correlations (joint probabilities) in the values of the underlying functions at any finite collection of points. Hence, this representation allows us to interpret probability measures on infinite-dimensional spaces as containing information about correlations across all possible finite collection of points in the spatial domain. Consequently, we can ``localize'' any infinite-dimensional probability measure. In particular, the first correlation marginal of this equivalent representation coincides with the classical notion of a Young measure. Thus, a probability measure on an $L^p$ space augments a Young measure with multi-point correlations and provides significantly more information than the Young measure does. We believe that this novel equivalence result could be of independent interest in stochastic analysis; see e.g.\ \cite{DZ}.

Another consequence of the equivalence of probability measures on function spaces and hierarchies of finite-dimensional correlation marginals, is the fact that the probability measure can be uniquely determined by a family of moments of the corresponding correlation marginals. Hence, the \emph{infinite-dimensional} Liouville or Hopf equation for statistical solutions, as proposed in \cite{FMRT}, can be replaced by an equivalent family of evolution equations (for moments) on \emph{finite-dimensional} (tensor-product) domains. 

The \emph{second aim} of this paper is to utilize this novel representation to define a suitable notion of statistical solutions for systems of conservation laws \eqref{eq:clcauchy}. In particular, certain moments (correlation functions) of the (time-parametrized) correlation marginals are evolved in a manner consistent with the dynamics of the conservation law \eqref{eq:clcauchy}. Consequently, statistical solutions need to satisfy an infinite family of evolutionary PDEs, but each of these PDEs is defined on a finite-dimensional spatial domain. 

The \emph{final aim} of this paper is to study the well-posedness of the proposed notion of statistical solutions. We will do so in the specific context of scalar conservation laws where we show existence of statistical solutions for a very large class of initial probability measures. The harder issue of uniqueness of statistical solutions for scalar conservation laws is also addressed. To this end, we propose a novel \emph{admissibility criterion} that amounts to requiring stability of each admissible statistical solution with respect to a specific set of stationary statistical solutions, namely those probability measures supported on finite collections of constant functions. Furthermore, we also show stability of the admissible statistical solution in the Wasserstein metric, with respect to probability measure-valued initial data: $W_1(\mu_t,\rho_t) \leq W_1(\bar{\mu},\bar{\rho})$. Thus, a complete characterization --- existence, uniqueness and stability --- of statistical solutions for scalar conservation laws is provided. The issues of existence and stability of admissible statistical solutions for the general case of systems of conservation laws will be presented in forthcoming papers in this series.

The rest of the paper is organized as follows. In Section \ref{sec:funcmeas} we prove the equivalence between probability measures on $L^p$ spaces and hierarchies of Young measures on finite-dimensional spaces. Statistical solutions for systems of conservation laws are defined in Section \ref{sec:stat} and the well-posedness of statistical solutions for scalar conservation laws is presented in Section \ref{sec:scalar}.

\section{Probability measures on function spaces}\label{sec:funcmeas}
The aim of this section is to establish the equivalence between probability measures on a function space, and families of measures describing the correlation of the values of underlying functions at different spatial points. The function spaces that we have in mind are $L^p(D,\U)$ for $1\leq p <\infty$ for some domain $D\subset\R^d$ and $\U:=\R^N$ (we will think of $D$ as physical space and $U$ as phase space). For ease of notation we will denote 
\[
\F:= L^p(D,\U).
\]
Henceforth, we equip $\F$ with its Borel $\sigma$-algebra $\Borel(\F)$.

A short summary of the contents this section follows. Given a probability measure $\mu$ on $\F = L^p(D,U)$, we might be interested in \emph{local} quantities such as the mean or the variance at a fixed point $x\in D$:
\[
\text{mean at }x = \int_\F u(x)\,d\mu(u), \qquad \text{variance at }x = \int_\F \bigl(u(x) - \text{mean})^2\,d\mu(u),
\]
or we might be interested in joint probability distributions at points $x, y\in D$:
\[
\text{probability that $u(x)\in A$ \textit{and} $u(y)\in B$} = \int_\F \ind_A(u(x))\ind_B(u(y))\,d\mu(u).
\]
However, not only are the integrands in the above integrals non-measurable, they are \textit{ill-defined} because point values $u(x)$ of a measurable function $u$ is not well-defined. Thus, we would like an equivalent representation of $\mu$ in terms of locally defined probability distributions $\nu^1_x$ or $\nu^2_{x,y}$; the above quantities could then be written as
\[
\int_U \xi\,d\nu^1_x(\xi), \qquad \int_U \Bigl(\xi-{\textstyle \int_U \xi'\ d\nu^1_x(\xi')}\Bigr)^2\,d\nu^1_x(\xi), \qquad \int_{U^2}\ind_A(\xi)\ind_B(\zeta)\,d\nu^2_{x,y}(\xi,\zeta)= \nu^2_{x,y}(A\times B),
\]
respectively. As we will see, we will require \textit{all} joint distributions across finitely many points in order to determine $\mu$ uniquely. This gives rise to an infinite hierarchy $\nu = (\nu^1,\nu^2,\dots)$ of maps $\nu^k$ from $D^k$ into $\Prob(U^k)$, the set of probability measures on $U^k$. Such a hierarchy is termed a \define{correlation measure} and each map $\nu^k$ a \define{correlation marginal}. The complete definition of correlation measures is given in Section \ref{sec:corrmeas}.

A similar construction is found in the Kolmogorov Extension Theorem (see e.g.\ \cite[Theorem 2.1.5]{Øks}). However, this approach considers measures on the product space $U^D$ equipped with the cylinder $\sigma$-algebra, instead of measures on $L^p(D,\U)$ equipped with its Borel $\sigma$-algebra. In the former case, questions such as ``is $u$ continuous'' or ``is $u$ Lebesgue integrable'' are not measurable, thus disqualifying its use in our context.

\subsection{Preliminaries}
We begin by recalling several definitions and results in functional analysis, measure theory and optimal transport theory.

\begin{notation}
If $\xi,\zeta\in\U$ then $\xi\cdot\zeta$ denotes their Euclidean inner product. If $D$ is a Borel set then $$D^k := \underbrace{D\times\dots\times D}_{k~{\rm times}}$$ and if $x=(x_1,\dots,x_k)\in D^k$ then we denote $|x| = |x_1|+\dots+|x_k|$. 

We denote the dual space of $\F$ by  $\F^* = L^{p'}(D,\U)$ (where $\frac{1}{p}+\frac{1}{p'}=1$), and the duality pairing between $\phi\in\F^*$ and $u\in\F$ by
\[
\ip{\phi}{u} = \phi(u) = \int_D \phi(x)\cdot u(x)\,dx.
\]

For any normed space $X$, we let $C_b(X)$ denote the space of bounded, continuous, real-valued functionals on $X$, equipped with the supremum norm $\|f\|_{C_b(X)} = \sup_{x\in X}|f(x)|$. We let $C_c(X)$ be the set of $f\in C_b(X)$ that have compact support, and we let $C_0(X)$ be the completion of $C_c(X)$ in the supremum norm.

The $k$-dimensional Lebesgue measure of a Borel set $A\subset\R^k$ is denoted $|A|$. The average of a function $f$ over a set $A$ is denoted
\[
\intavg_A f(x)\,dx = \frac{1}{|A|}\int_A f(x)\,dx.
\]

The Borel $\sigma$-algebra on a Polish space $X$ (i.e., a complete, separable metric space) is denoted by $\Borel(X)$. We let $\M(X)$ denote the space of finite, signed Radon measures on $(X,\, \Borel(X))$, and for $\mu\in\M(X)$ and $f\in L^1(X;\mu)$ we write $\ip{\mu}{f} = \int_X f(x)\,d\mu(x)$. The set $\Prob(X)$ of probability measures on $X$ consist of those $\mu\in\M(X)$ satisfying $\mu\geq0$ and $\mu(X)=1$.
\end{notation}

\subsubsection{The Wasserstein distance}
\begin{definition}\label{def:wasserstein}
Let $X$ be a separable Banach space and let $\mu, \rho\in\Prob(X)$ have finite $p$th moments, i.e.\ $\int_X |x|^p d\mu(x) < \infty$ and $\int_X |x|^p d\rho(x) < \infty$. The \define{$p$-Wasserstein distance} between $\mu$ and $\rho$ is defined as
\begin{equation}\label{eq:wasserdef}
W_p(\mu,\rho) = \inf_{\pi\in\Pi(\mu,\rho)}\int_{X^2} |x-y|^p\,d\pi(x,y);
\end{equation}
where the infimum is taken over the set $\Pi(\mu,\rho)\subset\Prob(X^2)$ of all transport plans from $\mu$ to $\rho$, i.e.\ those $\pi\in\Prob(X^2)$ satisfying
\[
\int_{X^2} F(x)+G(y)\,d\pi(x,y) = \int_X F(x)\,d\mu(x) + \int_X G(y)\,d\rho(y) \qquad \forall\ F,G\in C_b(X)
\]
(see e.g.\ \cite{Vil}). When $p=1$ we can write
\begin{equation}\label{eq:kantrubin}
W_1(\mu,\rho) = \sup_{\substack{\Psi\in C_b(X) \\ \|\Psi\|_{\Lip}\leq 1}} \int_X \Psi(x)\,d(\mu-\rho)(x),
\end{equation}
where the supremum is taken over all Lipschitz continuous functions with Lipschitz constant at most 1.
\end{definition}
It is straightforward to show that there always exists an \emph{optimal} transport plan $\pi$, i.e, one for which the infimum in \eqref{eq:wasserdef} is attained \cite[Theorem 1.3]{Vil}. The fact that \eqref{eq:wasserdef} and \eqref{eq:kantrubin} coincide when $p=1$ is a theorem in optimal transport theory often called the Kantorovich--Rubinstein theorem \cite[Theorem 1.14]{Vil}. The Wasserstein distance is a complete metric on the set of probability measures with finite $p$th moment, and metrizes the topology of weak convergence on this set \cite[Proposition 7.1.5]{AGS05}.

\subsubsection{Cylinder sets and -functions}
\begin{definition}
Let $X$ be a normed vector space. A function $\Psi : X \to \R$ is a \define{cylinder function} if there exist functionals $\phi_1,\dots,\phi_n \in X^*$ and a Borel measurable function $\psi : \R^n \to \R$ such that
\begin{equation}\label{eq:cylfunc}
\Psi(u) = \psi\bigl(\lip{\phi_1}{u}, \dots, \lip{\phi_n}{u}\bigr) \qquad \forall \ u\in X.
\end{equation}
A set $E\subset X$ is a \define{cylinder set} if the indicator function $u \mapsto \ind_E(u)$ is a cylinder function, or equivalently, if $E$ is of the form 
\begin{equation}\label{eq:cylsetdef}
E = \bigl\{ u\in X\ :\ \bigl(\lip{\phi_1}{u},\dots,\lip{\phi_n}{u}\bigr) \in F\bigr\}
\end{equation}
for a Borel set $F\subset\R^n$ and $\phi_1,\dots,\phi_n\in X^*$. We let $\Cyl(X)$ denote the collection of cylinder sets in $X$.
\end{definition}

\begin{proposition}\label{prop:uniquecyl}
%Let the normed vector space $X$ be either separable, or the dual of a separable space. 
Let $X$ be a separable normed vector space.
Then:
\begin{enumerate}[label=\it (\roman*)]
\item The $\sigma$-algebra generated by $\Cyl(X)$ is equal to $\Borel(X)$.
\item If $\mu$ is a (signed) measure on $(X,\Borel(X))$ such that $\mu(A) = 0$ for all cylinder sets $A$, then $\mu \equiv 0$.
\end{enumerate}
\end{proposition}
\begin{proof}
See the appendix.
\end{proof}

\subsection{Correlation measures}\label{sec:corrmeas}
\begin{definition}\label{def:corrmeas}
A \define{correlation measure} is a collection $\nu = (\nu^1, \nu^2, \dots)$ of maps $\nu^k : D^k \to \Prob(\U^k)$ satisfying the following properties:
\begin{enumerate}[label=\it (\roman*)]
\item\textit{Weak* measurability:} Each map $\nu^k : D^k \to \Prob(\U^k)$ is weak*-measurable, in the sense that the map $x \mapsto \ip{\nu^k_x}{f}$ from $x\in D^k$ into $\R$ is Borel measurable for all $f\in C_0(\U^k)$ and $k\in\N$. In other words, $\nu^k$ is a Young measure from $D^k$ to $\U^k$.
\item\textit{$L^p$-boundedness:} $\nu$ is $L^p$-bounded, in the sense that
\begin{equation}\label{eq:corrlpbound}
\int_{D} \ip{\nu^1_x}{|\xi|^p}\,dx < +\infty.
\end{equation}
\item\textit{Symmetry:} If $\sigma$ is a permutation of $\{1,\dots,k\}$ and $f\in C_0(\R^k)$ then $\ip{\nu^k_{\sigma(x)}}{f(\sigma(\xi))} = \ip{\nu^k_x}{f(\xi)}$ for a.e.\ $x\in D^k$. Here, we denote $\sigma(x) = \sigma(x_1,x_2,\ldots, x_k) = (x_{\sigma_1},x_{\sigma_2},\ldots,x_{\sigma_k})$. $\sigma(\xi)$ is denoted analogously.
\item\textit{Consistency:} If $f\in C_0(\U^k)$ is of the form $f(\xi_1,\dots,\xi_k) = g(\xi_1,\dots,\xi_{k-1})$ for some $g\in C_0(\U^{k-1})$, then $\ip{\nu^k_{x_1,\dots,x_k}}{f} = \ip{\nu^{k-1}_{x_1,\dots,x_{k-1}}}{g}$ for almost every $(x_1,\dots,x_k)\in D^k$.
\item\textit{Diagonal continuity (DC):} If $B_r(x) := \bigl\{y\in D\ :\ |x-y|<r\bigr\}$ then
\begin{equation}\label{eq:dcproperty}
\lim_{r\to0}\int_D\intavg_{B_r(x)}\ip{\nu^2_{x,y}}{|\xi_1-\xi_2|^p}\,dy\,dx = 0.
\end{equation}
\end{enumerate}
Each element $\nu^k$ is called a \define{correlation marginal}. We let $\Lp^{p} = \Lp^{p}(D,\U)$ denote the set of all correlation measures from $D$ to $\U$.
\end{definition}

\begin{remarks}\label{rem:corrmeas}~
\begin{enumerate}[label=\it (\roman*)]
\item By combining the properties of symmetry and consistency, the expected value with respect to $\nu^k_x$ of a function depending on $l<k$ parameters $\xi_{i_1}, \dots, \xi_{i_l}$, can be written in terms of $\nu^l_{x_{i_1}, \dots, x_{i_l}}$. Thus, the $k$th correlation marginal $\nu^k$ contains all information about lower-order correlation marginals, \emph{but not vice-versa}. Hence, the family $\nu = (\nu^k)_{k\in\N}$ constitutes a \emph{hierarchy}. 
\item Any function $u\in L^p(D,U)$ gives rise to a correlation marginal $\nu\in\Lp^p(D,U)$ by defining $\nu^k_x = \delta_{u(x_1)}\otimes\cdots\otimes\delta_{u(x_k)}$. Correlation marginals of this form are called \define{atomic}.
%\item The condition of diagonal continuity asserts (essentially) that $\nu^2$ equals $\nu^1$ along the diagonal $x=y$. By combining the condition of DC with the symmetry condition (and assuming that $\nu^k$ has been defined accordingly on the corresponding sets in $D^k$ of measure zero), we can write $\ip{\nu^k_x}{f}$ in terms of lower-order correlation marginals whenever one or more of the components of $x$ occur repeatedly.\alertUSF{Is this really needed?}
\item It can be shown that the DC property is equivalent to
\[
\lim_{r\to0}\int_D\intavg_{B_r(x)}\ip{\nu^2_{x,y}}{g(x,y,\xi_1,\xi_2)}\,dy\,dx = \int_D \ip{\nu^1_{x,x}}{g(x,x,\xi_1,\xi_1)}\,dx
\]
for every $g\in \Lone^2$. After possibly redefining $\nu^2$ on the zero-measure set $\{(x,y)\in D^2 : x=y\}$, this is equivalent to
\[
\nu^2_{x,x} = \nu^1_x \qquad \text{for a.e.\ }x\in D.
\]
%\item To illustrate the DC property, fix some arbitrary $\nu \in \Lp^2(D,U)$. Then a correlation measure would satisfy for all $f\in C_0(\U^2)$ and almost every $x\in D$,
%\[
%\ip{\nu^2_{x,x}}{f(\xi_1,\xi_2)} = \ip{\nu^1_x}{f(\xi,\xi)}.
%\]
In particular, $\ip{\nu^2_{x,x}}{\xi_1\xi_2} = \ip{\nu^1_x}{\xi^2}$ -- i.e., the covariance between the value at the point $x$ with itself is just the variance at $x$. Similarly, it can be shown that if $\psi\in C_b(U^{k+1})$ is Lipschitz continuous then
\[
\ip{\nu^{k+1}_{x_1,\dots,x_k,x_k}}{\psi(\xi_1,\dots,\xi_{k+1})} = \ip{\nu^{k}_{x_1,\dots,x_k}}{\psi(\xi_1,\dots,\xi_k,\xi_{k})}.
\]

We emphasize that diagonal continuity is an additional consistency requirement which is independent from consistency condition \textit{(iv)} of Definition \ref{def:corrmeas}.
\item As an example of a ``correlation measure'' which is not diagonally continuous, let $\nu^1 : D \to \Prob(\U)$ be any Young measure satisfying \eqref{eq:corrlpbound}, and define $\nu^k_{x_1,\dots,x_k} := \nu^1_{x_1}\otimes\cdots\otimes\nu^1_{x_k}$ for every $k\in\N$. Then $\nu = (\nu^1, \nu^2, \dots)$ satisfies properties \textit{(i)--(iv)} of Definition \ref{def:corrmeas}, but is DC \emph{if and only if} $\nu^1$ is atomic. Indeed, by Jensen's inequality,
\[
\ip{\nu^2_{x,x}}{\xi_1\xi_2} = \ip{\nu^1_x\otimes \nu^1_x}{\xi_1\xi_2} = \ip{\nu^1_{x}}{\xi}^2 \leq \ip{\nu^1_{x}}{\xi^2}
\]
for a.e.\ $x\in D$, with equality if and only if $\nu^1$ is atomic.
\end{enumerate}
\end{remarks}

\subsection{The main theorem}
Denote $\Lone^k := L^1\bigl(D^k,C_0(\U^k)\bigr)$. The proof of the following theorem, which is the main theorem of Section \ref{sec:funcmeas}, will depend crucially on $\Lone^k$ and its dual space; see Section \ref{sec:Lone}.
\begin{maintheorem}\label{thm:main}
For every correlation measure $\nu\in\Lp^p(D,\U)$ there exists a unique probability measure $\mu\in\Prob(\F)$ satisfying
\begin{equation}\label{eq:lpbound}
\int_\F \|u\|_\F^{p}\,d\mu(u) < \infty 
\end{equation}
such that
\begin{equation}\label{eq:nukdef}
%\int_{D^k}\ip{\nu^k_x}{g(x)}\,dx = \int_\F L_{g}(u)\,d\mu(u), \quad g\in\Cont^k.
\int_{D^k}\int_{\U^k}g(x,\xi)\,d\nu^k_x(\xi)dx = \int_\F \int_{D^k}g(x,u(x))\,dxd\mu(u) \qquad\forall\ g\in\Lone^k, \quad \forall\ k\in\N
\end{equation}
(where $u(x)$ denotes the vector $(u(x_1), \dots, u(x_k))$). Conversely, for every probability measure $\mu\in\Prob(\F)$ with finite moment \eqref{eq:lpbound}, there exists a unique correlation measure $\nu\in\Lp^p(D,\U)$ satisfying \eqref{eq:nukdef}. 

The relation \eqref{eq:nukdef} is also valid for any measurable $g:D\times\U\to\R$ such that $|g(x,\xi)|\leq C|\xi|^p$ for a.e.\ $x\in D$.
\end{maintheorem}
For a $g\in\Lone^k$, define the functional $L_g : \F\to\R$ by
\begin{equation}\label{eq:Lgdef}
L_g(u) := \int_{D^k}g(x,u(x))\,dx.
\end{equation}
Denoting $\ip{\nu^k}{g} := \int_{D^k}\int_{\U^k}g(x,\xi)\,d\nu^k_x(\xi)dx$, we can write \eqref{eq:nukdef} as
\begin{equation}\tag{\ref{eq:nukdef}'}\label{eq:nukdef2}
\ip{\nu^k}{g} = \ip{\mu}{L_g} \qquad \forall\ g\in\Lone^k, \quad \forall\ k\in\N.
\end{equation}
To ensure that the terms appearing in \eqref{eq:nukdef2} (or equivalently \eqref{eq:nukdef}) are well-defined, we need to check that $\nu^k$ is a continuous linear functional on $\Lone^k$, and that $L_g : \F\to\R$ is Borel measurable for every $g\in\Lone^k$. This is done in Theorem \ref{thm:Lonedual} and Proposition \ref{prop:Lgmeasurable}, respectively.

\begin{remark}
The finite moment requirement \eqref{eq:lpbound} is the direct analogue of the $L^p$ bound \eqref{eq:corrlpbound}. Indeed,
\[
\int_\F \|u\|_\F^{p}\,d\mu(u) = \int_\F \int_{D} |u(x)|^p\,dxd\mu(u) = \int_{D}\int_{\U} |\xi|^p\,d\nu^1_x(\xi) dx.\qedhere
\]
\end{remark}

\subsection{The spaces $\Lone^k$ and $\Lone^{k*}$}\label{sec:Lone}
\begin{definition}\label{def:Lone}
For any $k\in\N$, denote $\Lone^k := L^1\bigl(D^k,C_0(\U^k)\bigr)$, the space of measurable functions $g : x \mapsto g(x)\in C_0(\U^k)$ such that
\[
\|g\|_{\Lone^k} = \int_{D^k}\bigl\|g(x)\bigr\|_{C_0(\U^k)}\,dx < \infty.
\]
(Here, $C_0(\U^k)$ is equipped with its Borel $\sigma$-algebra.) We will routinely write $g(x,\xi)$ instead of $g(x)(\xi)$. We let $\Lone^{k*}:=L^\infty_w(D^k,\M(\U^k))$ denote the space of weak* measurable maps $\nu^k : x \mapsto \nu^k_x\in\M(\U^k)$ such that
\[
\|\nu^k\|_{\Lone^{k*}} = \esssup_{x\in D^k} \|\nu^k_x\|_{\M(\U^k)} < \infty.
\]
(Recall that $\nu^k$ is weak* measurable if the map $x \mapsto \ip{\nu^k_x}{f}$ from $D^k$ to $\R$ is measurable for all $f\in C_0(\U^k)$.)
\end{definition}

Note that if $\nu=(\nu^1,\nu^2,\dots)$ is a correlation measure then each correlation marginal $\nu^k$ is an element of $\Lone^{k*}$, because $\|\nu^k\|_{\Lone^{k*}} = 1$.

The following result justifies the notation $\Lone^{k*}$. 

\begin{theorem}\label{thm:Lonedual}
For any $k\in\N$, the space $\Lone^{k*}$ is isometrically isomorphic to the dual of $\Lone^k$ through the pairing
\[
\ip{\nu^k}{g} = \int_{D^k} \ip{\nu^k_x}{g(x,\cdot)}\,dx, \qquad g\in\Lone^k, \ \nu^k\in\Lone^{k*}.
\]
\end{theorem}
\begin{proof}
See e.g.\ \cite[Theorem 8.18.2]{Edw} or \cite[p.\ 211]{Bal89}.
\end{proof}

\begin{proposition}\label{prop:Lgmeasurable}
For any $g\in\Lone^k$, the map $L_g:\F\to\R$ defined by \eqref{eq:Lgdef} is uniformly continuous and satisfies
\begin{equation}\label{eq:Lgbounded}
\|L_g\|_{C_b(\F)} \leq \|g\|_{\Lone^k}.
\end{equation}
\end{proposition}
\begin{proof}
Since $g\in \Lone^k=L^1(D^k,C_0(\U^k))$, there are simple functions $\bar{g}_n(x) = \sum_{i=1}^{n}\ind_{A_{n,i}}(x)\bar{f}_{n,i}$ for functions $\bar{f}_{n,i}\in C_0(\U^k)$ and sets $A_{n,i}\subset D^k$ with positive and bounded Lebesgue measure, such that $\bar{g}_n \to g$ in $\Lone^k$. Let $f_{n,i}$ be functions in $C_0(\U^k)\cap \Lip(\U^k)$ such that $\|\bar{f}_{n,i}-f_{n,i}\|_{C_0(\U^k)} \leq \frac{1}{|A_{n,i}|n^2}$ (constructed, for instance, by mollification of $\bar{f}_{n,i}$), and define $g_n(x) := \sum_{i=1}^{n}\ind_{A_{n,i}}(x)f_{n,i}$.  If $u,v\in\F$ then
\begin{align*}
\bigl|L_{g_n}(u)-L_{g_n}(v)\bigr| &\leq \sum_{i=1}^{n}\int_{A_{n,i}}\bigl|f_{n,i}(u(x))-f_{n,i}(v(x))\bigr|\,dx \\
&\leq \sum_{i=1}^{n}\int_{A_{n,i}}\|f_{n,i}\|_{\Lip(\U^k)}\bigl(|u(x_1)-v(x_1)|+ \dots + |u(x_k)-v(x_k)|\bigr)\,dx \\
&\leq C_n\|u-v\|_\F
\end{align*}
by H\"older's inequality, where $C_n>0$ depends on $|A_{n,i}|$ and $\|f_{n,i}\|_{\Lip(\U^k)}$ for $i=1,\dots,n$. Thus, $L_{g_n}$ is Lipschitz continuous. Moreover,
\begin{align*}
|L_g(u)-L_{g_n}(u)| &\leq \int_{D^k}|g(x,u(x))-g_n(x,u(x))|\,dx \leq \int_{D^k} \|(g-g_n)(x)\|_{C_0(\U^k)}\,dx \\
&= \|g-g_n\|_{\Lone^k} \leq \|g-\bar{g}_n\|_{\Lone^k} + \frac{1}{n} \to 0 \qquad \text{as }n\to\infty,
\end{align*}
and so $L_{g_n} \to L_g$ uniformly on $\F$. Since every uniform limit of Lipschitz continuous functions is uniformly continuous, we conclude that $L_g$ is uniformly continuous. Finally,
\[
|L_g(u)| \leq \int_{D^k}|g(x,u(x))|\,dx \leq \int_{D^k}\|g(x)\|_{C_0(\U^k)}\,dx = \|g\|_{\Lone^k} \qquad \forall\ u\in\F,
\]
which proves \eqref{eq:Lgbounded}.
\end{proof}

\subsection{Existence and uniqueness of $\nu$}

\begin{theorem}\label{thm:nuexists}
Let $\mu\in\Prob(\F)$ satisfy \eqref{eq:lpbound}. Then \eqref{eq:nukdef} uniquely defines a correlation measure $\nu\in\Lp^p$.
\end{theorem}
\begin{proof}
We define each correlation marginal $\nu^k$ as an element of $\Lone^{k*}$ through duality, and then show that it has the required properties. The relation \eqref{eq:nukdef} uniquely defines $\nu^k$ as a linear functional on $\Lone^k$ which is continuous since
\[
|\ip{\nu^k}{g}| \leq \int_\F\int_{D^k}|g(x,u(x))|\,dxd\mu(u) \leq \int_{D^k}\|g(x)\|_{C_0(\U^k)}\,dx = \|g\|_{\Lone^k}.
\]
Thus, $\nu^k$ is an element of the dual of $\Lone^k$, which by Theorem \ref{thm:Lonedual} is $\Lone^{k*} := L^\infty_w(D^k,\M(\U^k))$. Hence, we can view $\nu^k$ as a weak* measurable map from $x\in D^k$ to $\nu^k_x \in \M(\U^k)$.

We show next that $\nu^k_x \in \Prob(\U^k)$ for Lebesgue-a.e.\ $x\in D^k$. For every $0\leq f\in C_0(\U^k)$ and for every bounded Borel measurable $A\subset D^k$ we have
\[
\ip{\nu^k}{\ind_A f} = \int_A \ip{\nu^k_x}{f}\,dx = \int_\F \int_A f(u(x_1),\dots,u(x_k))\,dxd\mu(u).
\]
But the right-hand side always lies between 0 and $|A|\cdot\|f\|_{C_0}$. It follows from the arbitrariness of $A$ that $0\leq \ip{\nu^k_x}{f}\leq \|f\|_{C_0}$ for Lebesgue-a.e.\ $x\in D^k$. In particular, letting $f(\xi) \equiv 1$, we find that $\|\nu^k_x\|_{\M} = 1$ for a.e.\ $x\in D$, which proves the claim.

Next, we show that $\nu=(\nu^1,\nu^2,\dots)$ satisfies properties \textit{(ii)--(iv)} of correlation measures (cf.\ Definition \ref{def:corrmeas}). The properties of symmetry and consistency follow directly from \eqref{eq:nukdef}, so it remains to show $L^p$-boundedness. By truncating the function $g:D\times U\to\R$ defined by $g(x,\xi) = |\xi_1|^p$ and applying Fatou's lemma and the dominated convergence theorem, we get that
\begin{align*}
\int_{D} \ip{\nu^1_x}{|\xi|^p}\,dx &= \ip{\nu^1}{|\xi|^p}
= \int_{\F}\int_{D} |u(x)|^p\,dxd\mu(u) = \int_{\F}\|u\|_{\F}^{p}\,d\mu(u) < +\infty.
\end{align*}
This proves \eqref{eq:corrlpbound}. 

Finally, we show that $\nu$ is diagonally continuous (cf.\ Definition \ref{def:corrmeas} \textit{(v)}). Indeed,
\begin{align*}
\lim_{r\to0}\int_D\intavg_{B_r(x)}\ip{\nu^2_{x,y}}{|\xi_1-\xi_2|^p}\,dy\,dx 
&= \lim_{r\to0}\int_\F\int_D\intavg_{B_r(x)}|u(x)-u(y)|^p\,dy\,dx\,d\mu(u) \\
&= \int_\F\int_D |u(x)-u(x)|^p\,dx\,d\mu(u) \\
&= 0,
\end{align*}
the second equality following from Lebesgue's differentiation theorem and the dominated convergence theorem. This completes the proof of existence of the correlation measure $\nu$. We emphasize that uniqueness follows directly from the explicit definition of $\nu^k$ (for each $k$) from \eqref{eq:nukdef}.
\end{proof}

\subsection{Uniqueness of $\mu$}

Let now $\nu\in\Lp^p(D,U)$ be a given correlation measure. We begin by proving that there exists at most one probability measure $\mu$ corresponding to $\nu$.

\begin{theorem} \label{lem:uniquefromcorr}
If $\mu,\tilde{\mu}\in\Prob(\F)$ both satisfy \eqref{eq:lpbound} and \eqref{eq:nukdef}, then $\mu = \tilde{\mu}$.
\end{theorem}
\begin{proof}
By assumption we have
\[
 \int_\F \int_{D^k}g(x,u(x))\,dxd\mu(u)
 = \int_\F \int_{D^k}g(x,u(x))\,dxd\tilde{\mu}(u)
  \qquad\forall\ g\in\Lone^k \quad \forall\ k\in\N.
\]
Fix a number $L>0$. By the dominated convergence theorem, H\"older's inequality and the $L^p$-bound \eqref{eq:lpbound}, this same equality holds for $g$ of the form 
$$
g(x,\xi) = \phi_1(x_1) \cdots \phi_k(x_k)\theta_L(x_1,\xi_1) \cdots \theta_L(x_k,\xi_k), \qquad \theta_L(x,\xi) = 
\begin{cases}
\xi & \text{if } |\xi| \leq L \text{ and } |x|\leq L \\
\frac{\xi}{|\xi|}L & \text{if } |\xi|>L \text{ and } |x|\leq L \\
0 & \text{if } |x|> L
\end{cases}
$$
where $\phi_1, \dots, \phi_k\in \F^*$. Denoting $\theta_L(u) = \theta_L(\cdot,u(\cdot))$ for the sake of simplicity, we can write \eqref{eq:nukdef} with the above test function $g$ as
\[
 \int_\F \langle \phi_1, \theta_L(u) \rangle \cdots \langle \phi_k, \theta_L(u) \rangle d\mu(u)
 = \int_\F  \langle \phi_1, \theta_L(u) \rangle \cdots \langle \phi_k, \theta_L(u) \rangle d\tilde{\mu}(u).
\]
By repeating indices (i.e.\ choosing some of the $\phi_i$'s to be identical) and expanding integrals over the spatial domain, one can show that the above identity implies
\begin{equation}\label{eq:moments}
 \int_\F \langle \phi_1, \theta_L(u) \rangle^{\alpha_1} \cdots \langle \phi_k, \theta_L(u) \rangle^{\alpha_k} d\mu(u)
 = \int_\F  \langle \phi_1, \theta_L(u) \rangle^{\alpha_1} \cdots \langle \phi_k, \theta_L(u) \rangle^{\alpha_k} d\tilde{\mu}(u).
\end{equation}
for arbitrary $\alpha_1, \dots, \alpha_k\in \N_0$.

Define now
\[
\phi:L^p(D) \to \R^k, \qquad 
\phi(u) := \Bigl(\ip{\phi_1}{u}, \dots, \ip{\phi_k}{u}\Bigr)
\]
and the truncation
$$
\phi_L:L^p(D) \to \R^k, \qquad 
\phi_L(u) := \Bigl(\ip{\phi_1}{\theta_L(u)}, \dots, \ip{\phi_k}{\theta_L(u)}\Bigr)
$$
Since $|\ip{\phi_i}{\theta_L(u)}| \leq m_d^{d/p}L^{1+d/p}\|\phi_i\|_{\F^*}$ for $i=1,\dots,k$ and with $m_d$ denoting the volume of the unit ball in $\R^d$, the map $\phi_L$ takes values only in the compact set $K_L := \bigl[-cL^{1+d/p},cL^{1+d/p}\bigr]^k\subset\R^k$, where $c = m_d^{d/p}\max\bigl(\|\phi_1\|_{\F^*},\dots,\|\phi_k\|_{\F^*}\bigr)$.

Let $\psi\in C^1_c(\R^k)$. Then the restriction of $\psi$ to $K_L$ can be approximated uniformly on $K_L$ by a sequence of polynomials $\bigl(P_n\bigr)_{n=1}^\infty$. It follows that 
$$
P_n\bigl(\phi_L(u)\bigr) \to \psi\bigl(\phi_L(u)\bigr) \qquad \text{as } n\to\infty
$$
uniformly in $u$. On the other hand, equation \eqref{eq:moments} implies that for each polynomial $P_n$, we have 
\[
\int_\F P_n\bigl(\phi_L(u)\bigr)\,d\mu(u) = \int_\F  P_n\bigl(\phi_L(u)\bigr)\,d\tilde{\mu}(u).
\]
From uniform convergence, we conclude that
\[
\int_\F \psi\bigl(\phi_L(u)\bigr)\,d\mu(u)
= \int_\F  \psi\bigl(\phi_L(u)\bigr)\,d\tilde{\mu}(u)
\]
for any $\psi \in C_c^1(\R^k)$. 

Define now $\Psi_L,\Psi : \F\to\R$ by
$$
\Psi_L(u) := \psi\bigl(\phi_L(u)\bigr), \qquad \Psi(u) := \psi\bigl(\phi(u)\bigr).
$$
Clearly, $|\Psi_L(u)|, |\Psi(u)| \leq \|\psi\|_{C_b(\R^k)}$ and $\lim_{L\to \infty} \Psi_L(u) = \Psi(u)$ for every $u\in\F$, so by the dominated convergence theorem,
\[
\int_\F \Psi(u) d\mu(u) = \int_\F  \Psi(u) d\tilde{\mu}(u)
\]
for any cylinder function $\Psi(u) = \psi\bigl(\ip{\phi_1}{u}, \dots, \ip{\phi_k}{u}\bigr)$ with $\psi\in C_c^1(\R^k)$. 

Given an open set $A \subset \R^k$, we can find a sequence $\psi_n \in C_c^1(\R^k)$ such that $0\leq\psi_n\leq\psi_{n+1}\leq\ind_A$ for all $n\in\N$, and $\psi_n$ converges pointwise to the indicator function $\ind_A$. Again, by dominated convergence, we conclude that 
\[
\int_\F \ind_A\bigl(\ip{\phi_1}{u}, \dots, \ip{\phi_k}{u}\bigr) d\mu(u)
 = \int_\F  \ind_A\bigl(\ip{\phi_1}{u}, \dots, \ip{\phi_k}{u}\bigr) d\tilde{\mu}(u).
\]
By a standard argument, this equality also holds for any Borel measurable set $A\subset\R^k$. This means that $\mu$ and $\tilde{\mu}$ agree on cylinder sets, so by Proposition \ref{prop:uniquecyl}, they must coincide.
\end{proof}

\subsection{Existence of $\mu$ for bounded $D$}\label{sec:muexistDbounded}
To prove existence of a probability measure $\mu$ corresponding to a given correlation measure $\nu$, we proceed in two steps, first proving the statement for bounded domains $D\subset \R^d$, and then extending the result to arbitrary $D\subset \R^d$.

We assume first that $D$ is bounded. Our construction will consist of a piecewise constant approximation over successively finer partitions of $D$.
\begin{definition}
A collection $\A = \{A_1,\dots,A_N\}$ of subsets of $D$ is a \define{partition of $D$} if
\[
\bigcup_{i=1}^N A_i = D,\qquad A_i\cap A_j = \emptyset \quad\text{and}\quad \bigl|\bar{A}_i\cap\bar{A}_j\bigr| = 0 \text{ for all } i\neq j
\]
(where $\bar{A}_i$ denotes the closure of $A_i$). Another partition $\widetilde{\A}=\big\{\widetilde{A}_1,\dots,\widetilde{A}_M\big\}$ is a \define{refinement of $\A$} if for every $j=1,\dots,M$, there is an $i\in\{1,\dots,N\}$ such that $\widetilde{A}_j \subset A_i$.
\end{definition}

Given a partition $\A = \{A_1, \dots, A_N\}$ of $D$ and a correlation measure $\nu\in\Lp^p(D,\U)$, define the probability measure $\rho_\A \in \Prob(\U^N)$ by
\[
\ip{\rho_\A}{\psi} = \intavg_{A_1\times\cdots\times A_N} \ip{\nu^N_x}{\psi}\,dx, \qquad \psi\in C_0(\U^N).
\]
This is clearly a continuous, linear functional on $C_0(\U^N)$ with norm $\displaystyle\|\rho_\A\|_{\M(\U^N)} = \sup_{\psi\in C_0}\frac{\ip{\rho_\A}{\psi}}{\|\psi\|_{C_0}} = 1$, and hence is a well-defined element of $\Prob(\U^N)$. Next, define $\mu_\A \in\Prob(\F)$ by
\[
\ip{\mu_\A}{\Psi} = \ip{\rho_\A}{\Psi\Bigl({\textstyle\sum_{i=1}^N\xi_i\ind_{A_i}}\Bigr)}.
\]
Being the pushforward of $\rho_\A$ by the continuous function $\U^N\ni \xi\mapsto\sum_{i=1}^N\xi_i\ind_{A_i}\in\F$, $\mu_\A$ is a well-defined element of $\Prob(\F)$. Finally, let $\nu_\A\in\Lp^p(D,\U)$ be the unique correlation measure corresponding to $\mu_\A$, as constructed in Theorem \ref{thm:nuexists}. It is clear that $\mu_\A$ is the probability measure corresponding to $\nu_\A$, in the sense of Theorem \ref{thm:main}. Note that $\nu_\A$ and $\mu_\A$ are piecewise constant, in the sense that each correlation marginal $\nu_{\A,x}^k$ is constant on sets of the form $x\in A_{i_1}\times\cdots\times A_{i_k}$, and $\mu_\A$ is concentrated on functions $u:D\to\U$ of the form $u(x)=\sum_{i=1}^N \xi_i \ind_{A_i}(x)$.

\begin{definition}
The correlation measure $\nu_\A\in\Lp^p(D,U)$ is called the \define{projection of $\nu$ onto $\A$}.
\end{definition}

It is not difficult to see that $\nu_\A$ can be equivalently defined as
\begin{equation}\label{eq:nuprojdef}
\ip{\nu^k_{\A,x}}{\psi} := \sum_{\alpha\in[N]^k} \ind_{A_\alpha}(x) \intavg_{A_1\times\cdots\times A_N} \ip{\nu^N_y}{\psi(\xi_{\alpha})}\,dy, \qquad x\in D^k \quad \forall\ k\in\N.
\end{equation}
(Here, $[N]=\{1,\dots,N\}$, $A_\alpha = A_{\alpha_1}\times\cdots\times A_{\alpha_k}$ and $\xi_\alpha = (\xi_{\alpha_1}, \dots,\xi_{\alpha_k})$.)

Given two partitions $\A$ and $\widetilde{\A}$ of $D$, where $\widetilde{\A}$ is a refinement of $\A$, the following lemma establishes an estimate for the distance between $\mu_{\A}$ and $\mu_{\widetilde{\A}}$.

\begin{lemma}\label{lem:proj_est}
Let $\nu\in\Lp^p(D,\U)$ be given. Let $\A$ and $\widetilde{\A}$ be partitions of $D$, where $\widetilde{\A}$ is a refinement of $\A$, and let $c,h>0$ be such that
\begin{equation}\label{eq:partitionassump}
|A_i| \geq ch^d, \qquad \mathrm{diam}(A_i) \leq h \qquad \forall\ A_i \in \A.
\end{equation}
Let $\mu_{\A}, \mu_{\widetilde{\A}}\in\Prob(\F)$ be the probability measures corresponding to the projections of $\nu$ onto $\A$ and $\widetilde{\A}$, respectively. Then
\[
W_1\bigl(\mu_\A,\mu_{\widetilde{\A}}\bigr) \leq C \Biggl(\int_{D}\intavg_{B_h(y)} \ip{\nu^2_{x,y}}{ |\xi_1-\xi_2|^p} \, dx\, dy\Biggr)^{1/p},
\]
where $B_h(y):=\{x\in D\ :\ |x-y|<h\}$ and $C>0$ only depends on $c$, $p$ and $d$ (the dimension of $D$).
\end{lemma}

\begin{proof}
Let $\Psi: \F \to \R$ be a Lipschitz function with $\|\Psi\|_{\Lip}=1$. Denote 
\[
\A = \big\{A_1, \dots, A_N\big\}, \qquad \widetilde{\A} = \big\{\widetilde{A}_1, \dots, \widetilde{A}_M\big\}.
\]
By definition,
\begin{align*}
\int_{\F} \Psi(u) \,d\bigl(\mu_{\A} - \mu_{\widetilde{\A}}\bigr) 
&= \intavg_{A_1\times\cdots\times A_N}\ip{\nu^N_{x}}{\Psi\Bigl({\textstyle\sum_{i=1}^N} \xi_i \ind_{A_i}\Bigr)}\,dx
- \intavg_{\widetilde{A}_1\times\cdots\times \widetilde{A}_M}\ip{\nu^M_y}{\Psi\Bigl({\textstyle\sum_{j=1}^M} \zeta_j \ind_{\widetilde{A}_j}\Bigr)}\,dy\\
&= \intavg_{A_1\times\cdots\times A_N}\ip{\nu^N_{x}}{\Psi\Bigl({\textstyle\sum_{j=1}^M} \xi_{i(j)} \ind_{A_{i(j)}}\Bigr)}\,dx
- \intavg_{\widetilde{A}_1\times\cdots\times \widetilde{A}_M}\ip{\nu^M_y}{\Psi\Bigl({\textstyle\sum_{j=1}^M} \zeta_j \ind_{\widetilde{A}_j}\Bigr)}\,dy\\
\end{align*}
where for any $j\in\{1, \dots,M\}$, the index $i(j)$ is the unique integer in $\{1,\dots, N\}$, such that $\widetilde{A}_{j}\subset A_{i(j)}$, and $\xi$ and $\zeta$ are the integration variables with respect to $\nu^N_x$ and $\nu^M_x$, respectively. Denote
\begin{gather*}
A = A_1\times \dots \times A_N, \qquad \widetilde{A} = \widetilde{A}_1 \times \dots \times \widetilde{A}_M.
\end{gather*}
Then we can write 
\begin{align*}
\int_{\F} \Psi(u) \, d\bigl(\mu_{\A} - \mu_{\widetilde{\A}}\bigr)
&= \intavg_{A} \ip{\nu^N_{x}}{\Psi \Bigl({\textstyle\sum_{j=1}^M} \xi_{i(j)} \ind_{\widetilde{A}_j}\Bigr)} \, dx  - \intavg_{\widetilde{A}} \ip{\nu^M_{y}}{\Psi \Bigl({\textstyle\sum_{j=1}^M} \zeta_j \ind_{\widetilde{A}_j}\Bigr)} \, dy \\
\text{\textit{(consistency of $\nu$)}}\quad &= \intavg_{\widetilde{A}}\intavg_{A} \ip{\nu^{N+M}_{x,y}}{\Psi \Bigl({\textstyle\sum_{j=1}^M} \xi_{i(j)} \ind_{\widetilde{A}_j}\Bigr)-\Psi \Bigl({\textstyle\sum_{j=1}^M} \zeta_j \ind_{\widetilde{A}_j}\Bigr)} \, dx\, dy \\
\text{\textit{(Lipschitz continuity)}}\quad &\leq \intavg_{\widetilde{A}}\intavg_{A} \ip{\nu^{N+M}_{x,y}}{\Bigl\| {\textstyle\sum_{j=1}^M} \xi_{i(j)} \ind_{\widetilde{A}_j}-{\textstyle\sum_{j=1}^M} \zeta_j \ind_{\widetilde{A}_j}\Bigr\|_{\F}} \, dx\, dy \\
&= \intavg_{\widetilde{A}}\intavg_{A} \ip{\nu^{N+M}_{x,y}}{\Bigl( {\textstyle\sum_{j=1}^M} |\widetilde{A}_j||\xi_{i(j)}-\zeta_j|^p \Bigr)^{1/p}} \, dx\, dy \\
\text{\textit{(Jensen's inequality)}}\quad &\leq \intavg_{\widetilde{A}}\intavg_{A} \left({\sum_{j=1}^M}|\widetilde{A}_j|\ip{\nu^{N+M}_{x,y}}{ |\xi_{i(j)}-\zeta_j|^p}\right)^{1/p} \, dx\, dy \\
\text{\textit{(Jensen's inequality)}}\quad &= \left(\intavg_{\widetilde{A}}\intavg_{A} {\sum_{j=1}^M}|\widetilde{A}_j|\ip{\nu^{N+M}_{x,y}}{ |\xi_{i(j)}-\zeta_j|^p} \, dx\, dy\right)^{1/p} \\
\text{\textit{(consistency of $\nu$)}}\quad &= \left(\sum_{j=1}^M|\widetilde{A}_j|\intavg_{\widetilde{A}}\intavg_{A} \ip{\nu^2_{x_{i(j)},y_j}}{ |\xi-\zeta|^p} \, dx\, dy\right)^{1/p} \\
&= \left(\sum_{j=1}^M|\widetilde{A}_j|\intavg_{\widetilde{A}_j}\intavg_{A_{i(j)}} \ip{\nu^2_{x_{i(j)},y_j}}{ |\xi-\zeta|^p} \, dx_{i(j)}\, dy_j\right)^{1/p}.
\end{align*}
Renaming variables $x_{i(j)}\mapsto x$ and $y_j\mapsto y$ in this summation, we obtain the estimate
\[
\int_{\F} \Psi(u) \, d\bigl(\mu_{\A} - \mu_{\widetilde{\A}}\bigr) \leq \left(\sum_{j=1}^M\int_{\widetilde{A}_j}\intavg_{A_{i(j)}} \ip{\nu^2_{x,y}}{ |\xi-\zeta|^p} \, dx\, dy\right)^{1/p},
\]
valid for any $1$-Lipschitz continuous $\Psi: L^1(D) \to \R$. Using \eqref{eq:partitionassump} we get the estimate 
\begin{align*}
\int_{\F} \Psi(u) \, d\bigl(\mu_{\A} - \mu_{\widetilde{\A}}\bigr)
&\leq \left(\sum_{j=1}^M\frac{1}{|A_{i(j)}|}\int_{\widetilde{A}_j}\int_{A_{i(j)}} \ip{\nu^2_{x,y}}{ |\xi-\zeta|^p} \, dx\, dy\right)^{1/p} \\
&\leq \left(\sum_{j=1}^M\frac{|B_h(y)|}{ch^d}\int_{\widetilde{A}_j}\intavg_{B_h(y)} \ip{\nu^2_{x,y}}{ |\xi-\zeta|^p} \, dx\, dy\right)^{1/p} \\
&\leq \left(C\int_{D}\intavg_{B_h(y)} \ip{\nu^2_{x,y}}{ |\xi-\zeta|^p} \, dx\, dy\right)^{1/p}
\end{align*}
where $C$ is given by the ratio of $c$ to the unit ball in $\R^d$. Taking the supremum over all $\Psi$ with $\|\Psi\|_{\Lip}\leq1$ on the left hand side and using the Kantorovich--Rubinstein definition \eqref{eq:kantrubin} of $W_1$ yields the desired estimate.
\end{proof}

With this bound in place we can complete the proof of existence of $\mu$.

\begin{theorem}
For any $\nu\in\Lp^p(D,\U)$ there exists a probability measure $\mu\in\Prob(\F)$ satisfying \eqref{eq:lpbound} and \eqref{eq:nukdef}.
\end{theorem}
\begin{proof}
Let $(\A_m)_{m\in\N}$ be a sequence of partitions of $D$ such that
\begin{itemize}
\item $\A_{m+1}$ is a refinement of $\A_m$,
\item there exists a constant $c>0$ and a sequence $h_m \to 0$, such that 
\[
|A|\geq ch_m^d, \qquad  \mathrm{diam}(A) \leq h_m \qquad \forall\ A\in \A_m, \; \forall\ m\in \mathbb{N}.
\]
\end{itemize}
We show first that the sequence of probability measures $\mu_{\A_m}\in\Prob(\F)$ converges weakly to some $\mu\in\Prob(\F)$ satisfying \eqref{eq:lpbound}. By Lemma \ref{lem:proj_est}, we have for any $m' > m$ 
\[
W_1\bigl(\mu_{\A_m},\mu_{\A_{m'}}\bigr) \leq C \Biggl(\int_{D}\intavg_{B_{h_m}(x)} \ip{\nu^2_{x,y}}{ |\xi-\zeta|^p} \, dy\, dx\Biggr)^{1/p}
\]
where $C>0$ does not depend on $m$. By the DC property \eqref{eq:dcproperty}, the right-hand side vanishes as $m\to\infty$. It follows that $\lim_{m,m'\to\infty}W_1\bigl(\mu_{\A_m},\mu_{\A_{m'}}\bigr) = 0$, so the sequence $\mu_{\A_m}$ is Cauchy in the $W_1$ metric. Since the $W_1$ metric turns $\Prob(\F)$ into a complete metric space (see \cite[Proposition 7.1.5]{AGS05}), we conclude that $\mu_{\A_m} \wto \mu$ for some $\mu\in\Prob(\F)$. Moreover, from the fact that $\nu$ satisfies \eqref{eq:corrlpbound}, it follows that $\mu$ satisfies \eqref{eq:lpbound}.

We show next that the limit $\mu$ satisfies \eqref{eq:nukdef}. Fix some $m\in\N$ and denote $\A = \A_m = \{A_1,\dots,A_N\}$. If $x\in D^k$ then there is a unique index $\alpha\in[N]^k$ such that $x\in A_\alpha:=A_{\alpha_1}\times\dots\times A_{\alpha_k}$. If $x$ is on the \emph{off-diagonal}, i.e. $\alpha_i\neq\alpha_j$ for all $i\neq j$, then it follows from consistency that
\[
\ip{\nu_{\A,x}^k}{\psi} = \sum_{\alpha\in[N]^k} \ind_{A_\alpha}(x) \intavg_{A_\alpha} \ip{\nu^k_y}{\psi}\,dy
\]
(compare with \eqref{eq:nuprojdef}). Hence, Lebesgue's differentiation theorem implies that $\ip{\nu_{\A_m,x}^k}{\psi} \to \ip{\nu_{x}^k}{\psi}$ as $m\to0$ for almost every point $x\in D^k$ on the off-diagonal $\bigl\{x\in D^k : x_i \neq x_j \text{ for all } i\neq j\bigr\}$. But since the diagonal $\bigl\{x\in D^k : x_i = x_j \text{ for some } i\neq j\bigr\}$ has Lebesgue measure zero, we can conclude that
\[
\wslim_{m\to\infty} \nu^k_{\A_m} = \nu^k \qquad \text{in } \Lone^{k*} \quad \forall\ k\in\N,
\]
or in other words,
\begin{equation}\label{eq:nuconv}
\lim_{m\to\infty} \ip{\nu^k_{\A_m}}{g} = \ip{\nu^k}{g} \qquad \forall\ g\in \Lone^k \quad \forall\ k\in\N.
\end{equation}

We know that $\mu_{\A_m} \wto \mu$ in $\Prob(\F)$, that is,
\begin{equation}\label{eq:muconv}
\lim_{m\to\infty}\int_\F \Psi(u)\,d\mu_{\A_m}(u) = \int_\F \Psi(u)\,d\mu(u) \qquad \forall\ \Psi\in C_b(\F).
\end{equation}
By Proposition \ref{prop:Lgmeasurable}, the functionals $L_g$ lie in $C_b(\F)$, so the above holds for $\Psi=L_g$ for any $g\in\Lone^k$. Thus, for any $k\in\N$ and $g\in\Lone^k$, we have
\[
\ip{\mu}{L_g} = \lim_{m\to\infty}\ip{\mu_{\A_m}}{L_g} = \lim_{m\to\infty}\ip{\nu^k_{\A_m}}{g} = \ip{\nu^k}{g},
\]
%We approximate now
%\begin{align*}
%\big|\ip{\mu}{L_g}-\ip{\nu^k}{g}\big| \leq \big|\ip{\mu-\mu_\A}{L_{g}}\big| +  \big|\ip{\mu_\A}{L_{g}}-\ip{\nu^k_\A}{g}\big| + \big|\ip{\nu^k_\A-\nu^k}{g}\big|
%\end{align*}
%The second term is zero, by definition of $\mu_\A$, and by \eqref{eq:nuconv} and \eqref{eq:muconv}, the first and third terms can be made small by letting $m$ be large enough. We conclude that
%\[
%\ip{\mu}{L_g} = \ip{\nu^k}{g} \qquad \forall\ g\in\Lone^k \quad \forall\ k\in\N,
%\]
which is \eqref{eq:nukdef}. 
%Finally, \eqref{eq:lpbound} follows from the dominated convergence theorem and \eqref{eq:nukdef}:
%\[
%\int_\F \|u\|_\F^p\,d\mu(u) = \int_D\ip{\nu^1_x}{|\xi|^p}\,dx < \infty.
%\]
\end{proof}

%\begin{remark}
%The above argument also implies that the limit probability measure $\mu$ constructed in this way does not depend on the particular choice of refinements $\A_m$. Indeed, from Lemma \ref{lem:proj_est}, we have that if each $\widetilde{\A}_m$ is a refinement of $\A_m$, then $W_1\bigl(\mu_{\widetilde{\A}_m}, \mu_{\A_m}\bigr) \to 0$. Given two such sequences, we can form a common refinement and conclude that they have the same limit.
%\end{remark}

\subsection{Existence of $\mu$ for unbounded $D$}

The next step is to prove existence of a probability measure $\mu$ for a given correlation measure $\nu$ on an arbitrary domain $D$. To this end, we first construct $\mu$ on a bounded set $E\subset D$, and then pass to the limit $E\uparrow D$.
\begin{lemma}\label{lem:restriction}
Let $E\subset D$. Let $r$ denote the restriction map
$$
r: L^p(D,\U) \to L^p(E,\U), \qquad  r(u) = u\bigr|_E.
$$
If $\mu\in \Prob\bigl(L^p(D,\U)\bigr)$ has correlation measure $\nu$, then $r\#\mu \in \Prob\bigl(L^p(E,\U)\bigr)$ has correlation measure 
\[
\nu\big|_E := \Bigl(\nu^1\big|_E,\, \nu^2\big|_{E^2},\, \nu^3\big|_{E^3},\, \dots\Bigr).
\]
\end{lemma}

\begin{proof}
Let $g\in L^1(E^k,C_0(U^k))$. Then the function $x \mapsto \ind_E(x)g(x,\cdot)$ lies in $\Lone^k=L^1(D^k,C_0(\U^k))$. Hence,
\begin{align*}
\int_{L^p(E,\U)} \int_{E^k} g(x,u(x)) \, dx \, d(r\#\mu)(u)
&= \int_{L^p(D,\U)} \int_{E^k}  g(x,u|_E(x)) \, dx \, d\mu(u)\\
&= \int_{L^p(D,\U)} \int_{D^k}  \ind_{E^k}(x) g(x,u(x)) \, dx \, d\mu(u) \\
&= \int_{D^k} \ip{\nu_x^k}{\ind_{E^k}(x)g(x,\cdot)} \, dx \\
&= \int_{E^k} \ip{\bigl(\nu^k\big|_{E^k}\bigr)_x}{g(x,\cdot)} \, dx.
\end{align*}
Thus, $\nu|_E$ is the correlation measure associated with $r\#\mu$.
\end{proof}

Let now $\nu\in\Lp^p(D,\U)$ for an arbitrary measurable set $D\subset\R^d$. Given $L>0$, let $D_L := D \cap (-L,L)^d$. Let $\tilde{\mu}_L\in\Prob(L^p(D_L,U))$ be the unique probability measure associated with the restriction $\nu|_{D_L}$ of $\nu$ to $D_L$, as constructed in Section \ref{sec:muexistDbounded}. Furthermore, let $\mu_L\in \Prob(L^p(D,\U))$ be the image of $\tilde{\mu}_L$ under the inclusion map obtained via extension by $0$:
$$
i_L: L^p(D_L,\U) \to L^p(D,\U),\qquad i_L(u) = u\ind_{D_L}.
$$
By Lemma \ref{lem:restriction}, we expect the sequence $(\tilde{\mu}_L)_{L>0}$ to be related to the restriction of a probability measure $\mu$ with correlation measure $\nu$. In particular, we would then expect the sequence $\mu_L$ to converge to a probability measure $\mu$ as $L\to\infty$. The following theorem shows that this is indeed the case.

\begin{theorem}
The sequence $\mu_L$ converges weakly as $L\to\infty$ to some $\mu\in\Prob(\F)$ satisfying \eqref{eq:lpbound} and \eqref{eq:nukdef}.
\end{theorem}
\begin{proof}
Let $\Psi\in C_b(\F)$ be an arbitrary $1$-Lipschitz function. Let $M < L,L'$. Then
\[
\int_{\F} \Psi(u) \, d\bigl(\mu_L - \mu_{L'}\bigr)
= \int_{\F} \Psi(u)-\Psi(\ind_{D_M}u) \, d\mu_L + \int_{\F} \Psi(\ind_{D_M}u) \, d\bigl(\mu_L - \mu_{L'}\bigr) + \int_{\F}\Psi(\ind_{D_M}u)-\Psi(u) \, d\mu_{L'}
\]
The second term is zero as a consequence of Lemma \ref{lem:restriction}. 
For the first and third terms, we have the estimate 
\begin{align*}
\left|\int_{\F} \Psi(u)-\Psi(\ind_{D_M}u) \, d\mu_L \right|
&\leq \int_{\F} \|u-\ind_{D_M}u\|_{L^p} \, d\mu_L \\
&\leq \Bigg(\int_{\F} \|\ind_{D_M^c}u\|_{L^p}^p \, d\mu_L\Bigg)^{1/p} \\
&= \Bigg(\int_{D_L \cap D_M^c} \langle \nu^1_x, | \xi |^p \rangle \, dx\Bigg)^{1/p} \\
&\leq \Bigg(\int_{D\setminus D_M} \langle \nu_x^1, | \xi |^p \rangle \, dx\Bigg)^{1/p}.
\end{align*}
It follows that
\[
\int_{\F} \Psi(u) \, d\bigl(\mu_L - \mu_{L'}\bigr)
\leq 2\left(\int_{D\setminus D_M} \langle \nu_x^1, | \xi |^p \rangle \, dx\right)^{1/p}
\]
Taking the supremum over all $1$-Lipschitz $\Psi\in C_b(\F)$ on the left, we obtain 
\[
W_1(\mu_L,\mu_{L'}) \leq 2\left(\int_{D\setminus D_M} \langle \nu_x^1, | \xi |^p \rangle \, dx\right)^{1/p}.
\]
By assumption, $\int_{D} \langle \nu_x^1, | \xi |^p \rangle \, dx$ is finite, so $\int_{D\setminus D_M} \langle \nu_x^1, | \xi |^p \rangle \, dx$ goes to zero as $M \to \infty$. We conclude that $W_1(\mu_L,\mu_{L'})\to 0$ as $L,L' \to \infty$. By completeness under the $1$-Wasserstein distance, the sequence $\mu_L$ converges to a limit $\mu = \wlim_{L \to \infty} \mu_L$.

We claim that the limit $\mu$ has correlation measure $\nu$, in the sense of Theorem \ref{thm:main}. Indeed, we have $\nu^k\bigl|_{D_L} \wsto \nu^k$, $\mu_L \wto \mu$ and $\ip{\nu^k\bigl|_{D_L}}{g} = \ip{\mu_L}{L_g}$ for all $g\in\Lone^k$ and $k\in\N$. It follows that $\ip{\nu^k}{g}=\ip{\mu}{L_g}$.
\end{proof}

\subsection{Moments}\label{sec:moments}
We have now established the equivalence between probability measures $\mu\in\Prob(\F)$ satisfying 
\begin{equation}\label{eq:moment1def}
\int_\F \|u\|_{\F}^p d\mu<\infty,
\end{equation}
and so-called correlation measures $\nu\in\Lp^p(D,U)$. In this section we introduce a third representation, that of moments. The \define{moments} of a correlation measure $\nu\in\Lp^p(D,U)$ are the functions
\begin{equation}\label{eq:momentdef}
m^k : D^k\to\U^{\otimes k}, \qquad m^k(x) := \int_{\U^k} \xi_1\otimes\cdots\otimes\xi_k\,d\nu^k_x(\xi), \qquad k\in\N.
\end{equation}
Here, $\U^{\otimes k}$ refers to the tensor product space $\U\otimes\cdots\otimes\U$ (repeated $k$ times), and $\xi_1\otimes\cdots\otimes\xi_k$ is a functional defined by its action on the dual space $\bigl(\U^{\otimes k}\bigr)^* = \U^{\otimes k}$ through
\[
\bigl(\xi_1\otimes\cdots\otimes\xi_k\bigr):\bigl(\zeta_1\otimes\cdots\otimes\zeta_k\bigr) = (\xi_1\cdot\zeta_1)\cdots(\xi_k\cdot\zeta_k).
\]
In the case $\U=\R$, the moments can be written more simply as
\[
m^k : D^k\to\R, \qquad m^k(x) = \int_{\R^k} \xi_1\cdots\xi_k\,d\nu^k_x(\xi), \qquad k\in\N.
\]
In either case, we will assume that
\begin{equation}\label{eq:nufinitemoments}
\int_{D^k} \int_{\U^k} |\xi_1|^p\cdots|\xi_k|^p\,d\nu^k_x(\xi) dx < \infty \qquad \forall\ k\in\N,
\end{equation}
or equivalently,
\begin{equation}\label{eq:mufinitemoments}\tag{\ref{eq:nufinitemoments}'}
\int_\F \|u\|_\F^{pk}\,d\mu(u) < \infty \qquad \forall\ k\in\N
\end{equation}
(compare with \eqref{eq:moment1def}). This ensures that $m^k$ is a well-defined element of $L^p(D^k,\, \U^{\otimes k})$.

The following result uniquely characterizes a correlation measure in terms of the family of moments $(m^k)_{k\in\N}$. This result will be essential to the contents of the following sections. 

\begin{theorem}\label{thm:momentunique}
Let $\nu\in\Lp^p(D,U)$ satisfy \eqref{eq:nufinitemoments}. Then the moments \eqref{eq:momentdef} uniquely identify $\nu$, in the sense that if another correlation measure $\tilde{\nu}$ has the same moments $(m^k)_{k\in\N}$, then $\nu=\tilde{\nu}$.
\end{theorem}
\begin{proof}
Denote by $\mu,\tilde{\mu}\in\Prob(\F)$ the corresponding probability measures. Recall that the \define{characteristic functional} of $\mu$ is the functional $\hat{\mu} : \F^*\to\R$,
\[
\hat{\mu}(\phi) := \int_\F e^{i\phi(u)}\,d\mu(u), \qquad \phi\in\F^*,
\]
and that $\mu$ and $\tilde{\mu}$ coincide if and only if $\hat{\mu}=\hat{\tilde{\mu}}$ (see \cite[Chapter 2.1]{DZ}). Using \eqref{eq:mufinitemoments} we can interchange integration and summation in the following and obtain
\begin{align*}
\hat{\mu}(\phi) &= \int_\F 1+\sum_{k=1}^\infty \frac{i^k}{k!}\phi(u)^k\,d\mu(u) = 1 + \sum_{k=1}^\infty \frac{i^k}{k!}\int_\F\left(\int_D\phi(x)\cdot u(x)\,dx\right)^kd\mu(u) \\
&= 1 + \sum_{k=1}^\infty \frac{i^k}{k!}\int_\F\int_{D^k}\bigl(\phi(x_1)\cdot u(x_1)\bigr)\cdots\bigl(\phi(x_k)\cdot u(x_k)\bigr)\,dxd\mu(u) \\
&= 1 + \sum_{k=1}^\infty \frac{i^k}{k!}\int_\F\int_{D^k}\bigl(u(x_1)\otimes\cdots\otimes u(x_k)\bigr) : \bigl(\phi(x_1)\otimes\cdots\otimes\phi(x_k)\bigr)\,dxd\mu(u) \\
&= 1 + \sum_{k=1}^\infty \frac{i^k}{k!}\int_{D^k}\int_{U^k}\bigl(\xi_1\otimes\cdots\otimes\xi_k\bigr):\bigl(\phi(x_1)\otimes\cdots\otimes\phi(x_k)\bigr)\,d\nu^k_xdx \\
&= 1 + \sum_{k=1}^\infty \frac{i^k}{k!}\int_{D^k}m^k(x):\bigl(\phi(x_1)\otimes\cdots\otimes\phi(x_k)\bigr)\,dx.
\end{align*}
Since the moments $m^k$ and $\tilde{m}^k$ of $\nu$ and $\tilde{\nu}$ coincide, we conclude that $\mu=\tilde{\mu}$.
\end{proof}

\subsection{Gaussian measures}
As an example of the equivalence of probability measures on function spaces and correlation measures, we present here a (somewhat formal) computation which characterizes the correlation measure for  \emph{Gaussian measures}, a class of probability measures that is of great interest in stochastic analysis \cite{DZ}. Although some of the following computations are quite standard in the literature on stochastic analysis, we include the details here for the sake of completeness.

We recall that a probability measure $\rho\in\Prob(\R)$ is \define{Gaussian} if there is a number $\sigma>0$ such that $\ip{\rho}{f} = \frac{1}{\sqrt{2\pi\sigma^2}}\int_\R f(z) e^{-\frac{z^2}{2\sigma^2}} dz$ for any $f\in C_0(\R)$. (Note that we are implicitly assuming that $\rho$ has mean zero, since the more general case of a nonzero mean can be easily obtained by translation.) Given a Banach space $X$, we say that a probability measure $\mu\in\Prob(X)$ is \define{Gaussian} if $\phi\#\mu \in \Prob(\R)$ is Gaussian for every nonzero $\phi\in X^*$, that is, if for every $0\neq \phi\in X^*$ there is a number $\sigma=\sigma(\phi)>0$ such that
\[
\int_X f(\phi(u))\,d\mu(u) = \frac{1}{\sqrt{2\pi\sigma^2}}\int_\R f(z)\exp\left(-\frac{z^2}{2\sigma^2}\right) dz \qquad \forall\ f\in C_0(\R).
\]
We easily find that the variance $\sigma(\phi)^2$ is given explicitly by
\begin{align*}
\sigma(\phi)^2 = \textrm{Var}(\phi\#\mu) &= \int_\R y^2\,d(\phi\#\mu)(y) = \int_X \phi(u)^2\,d\mu(u) = \ip{\mu}{\phi^2}.
\end{align*}

Choose now the Banach space $X=\F=L^p(D)$. For any $k\in\N$ and $0\neq\phi\in\F^*$, the expected value of the function $\R\ni z\mapsto z^k$ with respect to $\phi\#\mu$ is
\begin{align*}
\frac{1}{\sqrt{2\pi\sigma(\phi)^2}}\int_\R z^k \exp\left(-\frac{z^2}{2\sigma(\phi)^2}\right) dz &= \ip{\mu}{\phi^k}
= \int_\F \int_{D^k} \phi(x_1)u(x_1)\cdots\phi(x_k)u(x_k)\,dxd\mu(u) \\
&= \int_{D^k} \int_{\R^k} \xi_1\cdots\xi_k \phi(x_1)\cdots\phi(x_k)\,d\nu^k_x(\xi)dx \\
&= \int_{D^k} m^k(x)\phi(x_1)\cdots\phi(x_k)\,dx
\end{align*}
where $m^k(x):=\int_{\R^k}\xi_1\cdots\xi_k\,d\nu^k_x(\xi)$ denotes the $k$-th moment of $\nu$. On the other hand, it is well-known that the $k$-th moment $E[z^k]$ of a Gaussian distribution (with zero mean) is $0$ when $k$ is odd, and $(k-1)!!\sigma^k$ when $k$ is even, where $(k-1)!!$ denotes the double factorial $(k-1)!! = (k-1)(k-3)\cdots 1 = \frac{k!}{(k/2)!2^{k/2}}$. Using the fact that $m^k(x_1,\dots,x_k)$ is symmetric in all arguments, we find that $m^k\equiv0$ when $k$ is odd. When $k$ is even, i.e.\ $k = 2l$ for some $l\in\N$, we get
\begin{align*}
\int_{D^{2l}} m^{2l}(x)\phi(x_1)\cdots\phi(x_{2l})\,dx &= \frac{(2l)!}{l!2^l}\big(\sigma(\phi)^2\big)^l 
= \frac{(2l)!}{l!2^l}\left(\int_{D^{2}} m^2(x) \phi(x_1)\phi(x_2)\,dx\right)^l \\
&= \frac{(2l)!}{l!2^l}\int_{D^{2l}} m^2(x_1,x_2)\cdots m^2(x_{2l-1},x_{2l}) \phi(x_1)\cdots\phi(x_{2l})\,dx.
\end{align*}
The above implies that the first integrand must be given by the symmetric part of the last integrand, i.e.
\begin{align*}
m^{2l}(x) &= {\rm Sym}\left(\frac{(2l)!}{l!2^l}m^2\otimes\cdots\otimes m^2\right)(x) \\
&= \frac{1}{l!2^l}\sum_{s\in\mathfrak{S}_{2l}}m^2\bigl(x_{s(1)}, x_{s(2)}\bigr)\cdots m^2\bigl(x_{s(2l-1)}, x_{s(2l)}\bigr)
\end{align*}
where $\mathfrak{S}_{k}$ is the symmetric group on $k$ symbols, consisting of all permutations of $\{1, 2, \dots, k\}$ (see e.g.\ \cite{CGLM08}). Thus, all the moments---and thus all of $\mu$ (or, equivalently, $\nu$)---is completely specified in terms of the second moment $m^2$. (This general rule is known as Isserlis' theorem \cite{Iss18}; see also \cite[p.\ 44]{Frisch}.)

Finally, observe that 
\[
\ip{\nu^1_x}{\xi_1^n} = m^n(x,\dots,x)
= \begin{cases}
0 & \text{if $n$ is odd}, \\
(n-1)!!m^2(x,x)^{n/2} & \text{if $n$ is even}
\end{cases}
\]
(cf.\ Remark \ref{rem:corrmeas} \textit{(iii)}). Thus, for any $x\in D$, the probability measure $\nu^1_x$ is a Gaussian distribution with mean 0 and variance $m^2(x,x)$. More generally, for arbitrary $k$ we find that $\nu^k_{x_1,\dots,x_k}$ is a multivariate Gaussian distribution with mean $(0,\dots,0)$ and covariance $m^2(x_i,x_j)$. Thus, any function $m^2 : D^2 \to \R$ satisfying the properties of being a covariance function (see e.g.\ \cite[Section 21.2]{Kle}) corresponds to a unique Gaussian measure $\mu\in\Prob(L^p(D))$, and \textit{vice versa}. For instance, Brownian motion is obtained by letting $m^2(t,s)=\min(t,s)$ for $t,s\geq0$.

\section{Statistical solutions}
\label{sec:stat}
Equipped with the equivalence between probability measures on function spaces and correlation measures, we proceed in this section to define the concept of statistical solutions of multi-dimensional systems of conservation laws.

\subsection{Motivation and definition}
To motivate the equations governing the time-evolution of statistical solutions, we consider a scalar, one-dimensional conservation law
\[
\partial_t u + \partial_x f(u) = 0.
\]
This equation dictates the evolution of the quantity $u(x,t)$ over time. For $x_1,x_2\in \R$, consider the product $u(x_1,t)u(x_2,t)$. Assuming for the moment that $u$ is differentiable, we obtain
\begin{align*}
\partial_t\bigl[u(x_1,t)u(x_2,t)\bigr] &= \bigl(\partial_t u(x_1,t)\bigr)u(x_2,t) + u(x_1,t)\bigl(\partial_t u(x_2,t)\bigr) \\
&= -\partial_{x_1} f(u(x_1,t))u(x_2,t) - \partial_{x_2} u(x_1,t) f(u(x_2,t)),
\end{align*}
and for arbitrary $k\in\N$,
\begin{equation}\label{eq:momentpde}
\partial_t \bigl[u(x_1,t)\cdots u(x_k,t)\bigr] + \sum_{i=1}^k \partial_{x_i}\Bigl[u(x_1,t)\cdots f(u(x_i,t)) \cdots u(x_k,t)\Bigr] = 0.
\end{equation}
Since the above equation is in divergence form, it can be interpreted, in the sense of distributions, as
\begin{equation}\label{eq:weakcorrscalar}
\begin{split}
\int_{\R_+}\int_{\R^k} \partial_t \phi(x,t)\, u(x_1,t)\cdots u(x_k,t) + \sum_{i=1}^k \partial_{x_i}\phi(x,t)\, u(x_1,t)\cdots f(u(x_i,t)) \cdots u(x_k,t)\,dxdt \\
+ \int_{\R^k}\phi(0,x)\bar{u}(x_1)\cdots \bar{u}(x_k)\,dx = 0
\end{split}
\end{equation}
for all $\phi\in C_c^\infty(\R^k\times\R_+)$. 

For (multi-dimensional) systems, i.e.\ when $u$ and $f(u)$ are vectors, we evolve the tensor product $u(x_1)\otimes\cdots\otimes u(x_k)$, and the resulting evolution equation \eqref{eq:momentpde} would read
\begin{equation}\label{eq:weakcorrpde}
\partial_t \bigl[u(x_1,t)\otimes\cdots \otimes u(x_k,t)\bigr] + \sum_{i=1}^k \nabla_{x_i} \cdot \Bigl[u(x_1,t)\otimes\cdots \otimes f(u(x_i,t))\otimes \cdots\otimes u(x_k,t)\Bigr] = 0.
\end{equation}
Interpreting the above in the sense of distributions, we obtain
\begin{equation}\label{eq:weakcorr}
\begin{split}
\int_{\R_+}\int_{\R^k} \partial_t \phi(x,t): \bigl[u(x_1,t)\otimes\cdots \otimes u(x_k,t)\bigr] + \sum_{i=1}^k \nabla_{x_i}\cdot \phi(x,t) : \Bigl[u(x_1,t)\otimes\cdots \otimes f(u(x_i,t))\otimes \cdots\otimes u(x_k,t)\Bigr]\,dxdt \\
+ \int_{\R^k}\phi(0,x):\bigl[\bar{u}(x_1)\otimes\cdots\otimes\bar{u}(x_k)\bigr]\,dx = 0
\end{split}
\end{equation}
for all $\phi\in C_c^\infty\Bigl(\big(\R^d\big)^k\times\R_+,\ \big(\R^N\big)^{\otimes k}\Bigr)$.
The above calculations can be made rigorous, as follows.
\begin{lemma}\label{lem:weaksoln}
If $u\in L^1_{{\rm loc}}(\R^d\times\R_+,\,\R^N)$ is a weak solution of \eqref{eq:clcauchy} then \eqref{eq:weakcorr} holds for all $k\in\N$.
\end{lemma}
\begin{proof}
For the sake of notational simplicity we present the proof only for the one-dimensional, scalar case ($d=N=1$). The proof proceeds by induction. Equation \eqref{eq:weakcorr} with $k=1$ is precisely the definition of a weak solution,
\begin{equation}\label{eq:weaksoln}
\int_{\R_+}\int_{\R} \partial_t \psi u + \partial_x\psi u\,dxdt + \int_{\R}\psi(x,0)u(x)\,dx = 0 \qquad \forall\ \psi\in C_c^\infty(\R\times\R_+).
\end{equation}
Assume that \eqref{eq:weakcorr} holds for some $k\in\N$. Let $\omega_\eps:\R\to\R$ be a symmetric mollifier with $\supp\omega_\eps\subset[-\eps,\eps]$, let $\tilde{\phi}\in C_c^\infty(\R^{k+1}\times\R_+)$ and define
\[
\phi(x,t) := \int_{\R_+}\int_\R \omega_\eps(t-s)\tilde{\phi}(x,x_{k+1},s)u(x_{k+1},s)\,dx_{k+1}ds
\]
for $x\in \R^k$ and any $0\leq\tilde{\phi}\in C_c^\infty(\R^{k+1}\times\R_+)$. Then $\phi\in C_c^\infty(\R^k\times\R_+)$, and we have
\begin{align*}
\partial_t\phi(x,t) 
&= \int_{\R_+}\int_\R \omega_\eps'(t-s)\tilde{\phi}(x,x_{k+1},s)u(x_{k+1},s)\,dx_{k+1}ds \\
&= \int_{\R_+}\int_\R \left[-\partial_s\Bigl(\omega_\eps(t-s)\tilde{\phi}(x,x_{k+1},s)\Bigr) + \omega_\eps(t-s)\partial_s\tilde{\phi}(x,x_{k+1},s)\right]u(x_{k+1},s)\,dx_{k+1}ds \\
&= \int_{\R_+}\int_\R \omega_\eps(t-s)\partial_{x_{k+1}}\tilde{\phi}(x,x_{k+1},s)f(u(x_{k+1},s))\,dx_{k+1}ds \\
&\quad + \int_\R \omega_\eps(t)\tilde{\phi}(x,x_{k+1},0) \bar{u}(x_{k+1}) \,dx_{k+1}\\
&\quad + \int_{\R_+}\int_\R \omega_\eps(t-s)\partial_s\tilde{\phi}(x,x_{k+1},s)u(x_{k+1},s)\,dx_{k+1}ds,
\end{align*}
the last equality following from \eqref{eq:weaksoln}. Moreover, for $j=1,\dots,k$ we have
\[
\partial_{x_j}\phi(x,t) = \int_{\R_+}\int_\R \omega_\eps(t-s)\partial_{x_j}\tilde{\phi}(x,x_{k+1},s)u(x_{k+1},s)\,dx_{k+1}ds.
\]
Hence, inserting $\phi$ into \eqref{eq:weakcorr} gives
\begin{align*}
0 &= \int_{\R_+}\int_{\R^k} u(x_1,t)\cdots u(x_k,t)\Biggl[\int_{\R_+}\int_\R \omega_\eps(t-s)\partial_{x_{k+1}}\tilde{\phi}(x,x_{k+1},s)f(u(x_{k+1},s))\,dx_{k+1}ds \\
&\quad + \int_\R \omega_\eps(t)\tilde{\phi}(x,x_{k+1},0) \bar{u}(x_{k+1}) \,dx_{k+1} + \int_{\R_+}\int_\R \omega_\eps(t-s)\partial_s\tilde{\phi}(x,x_{k+1},s)u(x_{k+1},s)\,dx_{k+1}ds\Biggr] \\
&\quad+\sum_{j=1}^{k} u(x_1,t)\cdots f(u(x_j,t))\cdots u(x_k,t)\int_{\R_+}\int_\R \omega_\eps(t-s)\partial_{x_j}\tilde{\phi}(x,x_{k+1},s)u(x_{k+1},s)\,dx_{k+1}dsdxdt \\
&\quad + \int_{\R^k} \bar{u}(x_1)\cdots\bar{u}(x_k)\int_{\R_+}\int_\R \omega_\eps(-s)\tilde{\phi}(x,x_{k+1},s)u(x_{k+1},s)\,dx_{k+1}dsdx.
\end{align*}
In the limit $\eps\to0$ we get
\begin{align*}
0 &= \int_{\R_+}\int_{\R^k}\int_\R u(x_1,t)\cdots u(x_k,t) \partial_{x_{k+1}}\tilde{\phi}(x,x_{k+1},t)f(u(x_{k+1},t))\,dx_{k+1}dxdt \\
&\quad + \frac{1}{2}\int_{\R^k}\int_\R\bar{u}(x_1)\cdots\bar{u}(x_k) \tilde{\phi}(x,x_{k+1},0) \bar{u}(x_{k+1}) \,dx_{k+1}dx \\
&\quad + \int_{\R_+}\int_{\R^k}\int_\R u(x_1,t)\cdots u(x_k,t) \partial_t\tilde{\phi}(x,x_{k+1},t)u(x_{k+1},t)\,dx_{k+1}dxdt \\
&\quad+\int_{\R_+}\int_{\R^k}\int_\R\sum_{j=1}^{k} u(x_1,t)\cdots f(u(x_j,t))\cdots u(x_k,t) \partial_{x_j}\tilde{\phi}(x,x_{k+1},t)u(x_{k+1},t)\,dx_{k+1}dxdt \\
&\quad + \frac{1}{2}\int_{\R^k}\int_\R \bar{u}(x_1)\cdots\bar{u}(x_k) \tilde{\phi}(x,x_{k+1},0)\bar{u}(x_{k+1})\,dx_{k+1}dx
\end{align*}
(The factors $\frac{1}{2}$ come from integrating $\omega_\eps(-s)$ over $s\in\R_+$ and not $s\in\R$.) After reorganizing terms, we obtain \eqref{eq:weakcorr} for $k+1$.
\end{proof}

Denoting the atomic correlation measure corresponding to $u(\cdot,t)$ by $\nu_t = (\nu^{1}_t, \nu^{2}_t, \dots)$ (cf.\ Remark \ref{rem:corrmeas}\textit{(ii)}), we may write \eqref{eq:weakcorrpde} equivalently as
\begin{equation}\label{eq:momentcorrmeas}
\partial_t \ip{\nu^k_{t,x}}{\xi_1\otimes\cdots\otimes\xi_k} + \sum_{i=1}^k \nabla_{x_i}\cdot\ip{\nu^k_{t,x}}{\xi_1\otimes\cdots\otimes f(\xi_i)\otimes\cdots\otimes\xi_k} = 0
\end{equation}
for $x\in \R^k$, $t>0$ and any $k\in\N$. Note that this expression makes sense even if $\nu^k_t$ is non-atomic. We take this as the definition of a possibly non-atomic statistical solution. In order for the terms appearing in \eqref{eq:momentcorrmeas} to be well-defined, we need to assume
\[
\int_{\bar{D}^k} \ip{\nu^k_{t,x}}{|\xi_1|\cdots|\xi_k|}\,dx < \infty, \qquad \int_{\bar{D}^k} \ip{\nu^k_{t,x}}{|\xi_1|\cdots|f(\xi_i)|\cdots |\xi_k|}\,dx < \infty \qquad \forall\ k\in \N \quad {\rm and} \quad i = 1,2,\ldots,k,
\]
for all compact subsets $\bar{D} \subset D$.

We can write this in terms of the corresponding probability measure $\mu_t\in\Prob(L^1)$ as
\begin{equation}\label{eq:fastdecay}
\int_{L^1} \|u\|_{L^1(\bar{D}^k)}\,d\mu_t(u) < \infty, \qquad \int_{L^1(\bar{D}^k)} \|f\circ u\|_{L^1(\bar{D}^k)}\|u\|_{L^1(\bar{D}^k)}^k d\mu_t(u) < \infty \qquad \forall\ k\in\N,
\end{equation}
and for all compact subsets $\bar{D} \subset D$.
\begin{definition}\label{def:statsoln}
Let $\bar{\mu}\in\Prob\big(L^1\big(\R^d,\R^N\big)\big)$ satisfy the decay rate \eqref{eq:fastdecay}. A \define{statistical solution} of \eqref{eq:cl} with initial data $\bar{\mu}$ is a weak*-measurable mapping $t \mapsto \mu_t\in \Prob\big(L^1\big(\R^d,\R^N\big)\big)$ such that each $\mu_t$ satisfies the decay rate \eqref{eq:fastdecay}, and such that the corresponding correlation measures $(\nu^k_t)_{k\in\N}$ satisfy \eqref{eq:momentcorrmeas} in the sense of distributions, i.e.
\begin{gather*}
\int_{\R_+}\int_{(\R^d)^k} \ip{\nu^k_{t,x}}{\xi_1\otimes\cdots\otimes\xi_k}:\partial_t \phi + \sum_{i=1}^k \ip{\nu^k_{t,x}}{\xi_1\otimes\cdots\otimes f(\xi_i)\otimes\cdots\otimes\xi_k}: \nabla_{x_i}\phi\,dxdt \\
+ \int_{(\R^d)^k}\ip{\bar{\nu}^k_x}{\xi_1\otimes\cdots\otimes\xi_k}\phi\bigr|_{t=0}\,dx = 0
\end{gather*}
for every $\phi\in C_c^\infty\Big(\big(\R^d\big)^k\times\R_+,\ \big(\R^N\big)^{\otimes k}\Big)$ and for every $k\in\N$. We denote $\bar{\nu}$ to the correlation measure associated with initial probability measure $\bar{\mu}$.
\end{definition}
(A map $\mu : t\mapsto \mu_t\in\Prob\big(L^1\big(\R^d,\R^N\big)\big)$ is \define{weak*-measurable} if the pairing $\ip{\mu_t}{G} = \int_{L^1}G(u)\,d\mu_t(u)$ with any $G\in C_b\big(L^1\big(\R^d,\R^N\big)\big)$ is Lebesgue measurable in $t$ (see e.g.\ \cite[Section II.1]{DU77}).)

\begin{remarks}~
\begin{enumerate}[label=\it (\roman*)]
\item
Note carefully that the evolution equation \eqref{eq:momentcorrmeas} dictates the evolution of the \emph{moments} $\ip{\nu^k_{t,x}}{\xi_1\otimes\cdots\otimes\xi_k}$ (see Section \ref{sec:moments}). Recall from Theorem \ref{thm:momentunique} that the moments of a correlation measure uniquely identify the correlation measure. Thus, instead of determining the time evolution of functionals on infinite-dimensional function spaces as in the Liouville and Hopf equations of \cite{FMRT}, we reduce the problem to the evolution of functions $\ip{\nu^k_{t,x}}{\xi_1\otimes\cdots\otimes\xi_k}$ defined on the finite-dimensional spaces $(x,t)\in \big(\R^d\big)^k\times\R_+$.
\item
Equation \eqref{eq:momentcorrmeas} for $k=1$ is simply the definition of $\nu^{1}$ being a measure-valued solution of \eqref{eq:cl}, as introduced by DiPerna \cite{DiP85}. In light of the previous remark, we see that---except when the correlation measure is atomic---the evolution equation for measure-valued solutions (i.e., \eqref{eq:momentcorrmeas} with $k=1$) \emph{never uniquely determines} the full correlation measure $\nu_t$ (or equivalently, $\mu_t$). In other words, except in the case of an atomic statistical solution, the evolution equation for the $(k+1)$th moment can contain strictly more information than the equation for the $k$th moment. Thus, statistical solutions are much more constrained than measure-valued solutions with additional information being provided by multi-point correlation measures. This additional information provided by the correlation measures, opens the possibility of enforcing uniqueness of the statistical solutions, if necessary by augmenting them with further admissibility conditions. 
\item If $\bar{\mu} = \delta_{\bar{u}}$ and $\mu_t = \delta_{u(t)}$ with $\bar{u}, u(t) \in L^1\big(\R^d,\R^N\big)$ for a.e.\ $t>0$, then Definition \ref{def:statsoln} reduces to the classical definition of a weak solution of \eqref{eq:cl}.
\end{enumerate}
\end{remarks}

\section{Statistical solutions for scalar conservation laws}\label{sec:scalar}
In Section \ref{sec:stat} we defined statistical solutions for multi-dimensional \textit{systems} of conservation laws. In this section we investigate the well-posedness of statistical solutions of (multi-dimensional) \textit{scalar} conservation laws. To this end, we can utilize the well-posedness of the deterministic problem \eqref{eq:clcauchy} to show existence of a statistical solution for a multi-dimensional scalar conservation law.

\subsection{The canonical statistical solution}\label{sec:canonss}
Recall that for scalar conservation laws, the Cauchy problem \eqref{eq:clcauchy} is well-posed for any $\bar{u}\in\EuScript{U}:= L^1\cap L^\infty(\R^d\times\R_+)$, and the entropy solution $u(t) = S_t\bar{u}$ lies in $\EuScript{U}$ for all $t>0$ \cite{Kruz1}. Here, $S_t : \EuScript{U}\to\EuScript{U}$ denotes the entropy solution semi-group. Denote $\F := L^1(\R^d)$. Given initial data $\bar{\mu}\in\Prob(\F)$ with $\supp\mu\subset\EuScript{U}$, we define the \define{canonical statistical solution} by
\[
\mu_t := S_t\#\bar{\mu}, \qquad t\geq0,
\]
where the pushforward operator $\#$ applies $S_t$ to each element of the support of $\bar{\mu}$:
\[
\int_\F G(u)\,d\left(S_t\#\bar{\mu}\right)(u) = \int_\F G(S_tu)\,d\bar{\mu}(u), \qquad G\in C_b(\F).
\]
Thus, the canonical statistical solution is concentrated on the entropy solutions of every initial data in the support of $\bar{\mu}$, and each entropy solution is given the same weight as $\bar{\mu}$ gives to the corresponding initial data. 

The semi-group $S_t$ is a continuous map, so it is easy to see that the canonical statistical solution is a weak*-measurable map from $t\in\R_+$ to $\Prob(\F)$. Moreover, it is in fact a statistical solution: For every $k\in\N$ and $\phi\in C_c(\R^k\times\R_+)$, we have
\begin{align*}
&\int_{\R_+}\int_{(\R^d)^k} \partial_t \phi\, \ip{\nu^k_{t,x}}{\xi_1\cdots\xi_k} + \sum_{i=1}^k \nabla_{x_i}\phi\,: \ip{\nu^k_{t,x}}{\xi_1\cdots f(\xi_i)\cdots\xi_k}\,dxdt + \int_{(\R^d)^k}\phi\bigr|_{t=0}\,\ip{\bar{\nu}^k_x}{\xi_1\cdots\xi_k}\,dx \\
=&\ \int_{\R_+}\int_\F\int_{(\R^d)^k} \partial_t \phi\, u(x_1)\cdots u(x_k) + \sum_{i=1}^k \nabla_{x_i}\phi : \Bigl[u(x_1)\cdots f\bigl(u(x_i)\bigr) \cdots u(x_k)\Bigr]\,dxd\mu_t(u)dt \\
&+ \int_\F\int_{(\R^d)^k}\phi(0,x)\bar{u}(x_1)\cdots \bar{u}(x_k)\,dxd\bar{\mu}(\bar{u}) \qquad\qquad ({\rm by ~\eqref{eq:nukdef}}) \\
=&\ \int_\F \Biggl[\int_{\R_+}\int_{(\R^d)^k} \partial_t \phi\, S_t\bar{u}(x_1)\cdots S_t\bar{u}(x_k) + \sum_{i=1}^k \nabla_{x_i} \phi : \Bigl[S_t\bar{u}(x_1)\cdots f\bigl(S_t\bar{u}(x_i)\bigr) \cdots S_t\bar{u}(x_k)\Bigr]\,dxdt \\
&+ \int_{(\R^d)^k}\phi(0,x)\bar{u}(x_1)\cdots \bar{u}(x_k)\,dx\Biggr]d\bar{\mu}(\bar{u}) \\
=&\ 0
\end{align*}
by Lemma \ref{lem:weaksoln}, since $S_t \bar{u}$ is a weak solution of \eqref{eq:clcauchy} for every $\bar{u}\in\EuScript{U}$.

It is also quite easy to see that the canonical statistical solution is stable with respect to the initial data. We measure this stability in the 1-Wasserstein metric $W_1$ on $\F$ (cf.\ Definition \ref{def:wasserstein}). Let $\bar{\mu}, \bar{\rho} \in \Prob(\F)$ be given initial data and let $\bar{\pi}\in\Pi(\bar{\mu},\bar{\rho})$ be an optimal transport plan from $\bar{\mu}$ to $\bar{\rho}$. For each $t\geq0$ we define $\pi_t := (S_t,S_t)\#\bar{\pi}$, which lies in $\Pi(\mu_t,\rho_t)$ (where $\mu_t,\rho_t$ are the corresponding canonical statistical solutions). We find that
\begin{align*}
W_1(\mu_t,\rho_t) &\leq \int_{\F^2} \|u-v\|_{\F}\,d\pi_t(u,v) = \int_{\F^2} \|S_t\bar{u}-S_t\bar{v}\|_{\F}\,d\bar{\pi}(\bar{u},\bar{v}) \\
&\leq \int_{\F^2} \|\bar{u}-\bar{v}\|_{\F}\,d\bar{\pi}(\bar{u},\bar{v}) = W_1(\bar{\mu},\bar{\rho}),
\end{align*}
where the first inequality comes from picking a particular plan $\pi_t\in\Pi(\mu_t,\rho_t)$ in \eqref{eq:wasserdef}, and the second inequality follows from the $L^1$ contraction property of $S_t$. We summarize these observations as follows.
\begin{theorem}
Let $\bar{\mu}\in\Prob(\F)$ be a probability measure on $\F$ satisfying \eqref{eq:fastdecay}, and define the \define{canonical statistical solution} $\mu_t := S_t\#\bar{\mu}$ for each $t\in\R_+$. Then $t\mapsto \mu_t$ is a statistical solution of \eqref{eq:cl} with data $\bar{\mu}$, and if $\rho_t$ is another canonical statistical solution with initial data $\bar{\rho}\in\Prob(\F)$ then
\begin{equation}
W_1(\mu_t,\rho_t) \leq W_1(\bar{\mu},\bar{\rho}).
\end{equation}
\end{theorem}

\subsection{Well-posedness of statistical solutions}
As shown in Section \ref{sec:canonss}, there always exists a statistical solution for scalar conservation laws, and this solution is stable with respect to initial data. This does not imply, however, that the canonical solution is unique, in the same way that there might exist several weak solutions for the deterministic equation \eqref{eq:clcauchy}. As in the deterministic setting, \emph{entropy conditions} must be imposed in order to single out a unique solution.

Recall that the (Kruzkov) entropy condition for \eqref{eq:clcauchy} is
\begin{equation}\label{eq:kruzkovec}
\partial_t |u-c| + \nabla_x \cdot q(u,c) \leq 0 \qquad \text{in } \mathcal{D}'(\R^d\times\R_+)
\end{equation}
for all constants $c\in\R$, where $q(u,c) := \sgn(u-c)(f(u)-f(c))$. Although not usually phrased as such, the Kruzkov entropy condition imposes stability with respect to a certain family of stationary (steady-state) solutions, namely the constant solutions. The key to proving uniqueness of statistical solutions lies in finding the right family of stationary (time-invariant) solutions. A natural first attempt follows from integrating \eqref{eq:kruzkovec} over the phase-space variable, which yields
\begin{equation}\label{eq:entrcond1}
\partial_t \ip{\nu^1}{|\xi-c|} + \nabla_x \cdot \ip{\nu^1}{q(\xi,c)} \leq 0 \qquad \text{in } \mathcal{D}'(\R^d\times\R_+).
\end{equation}
This is the entropy condition enforced by DiPerna in the context of measure-valued solutions \cite{DiP85}. By a standard doubling-of-variables argument (see \cite[Theorem 4.1]{DiP85} and \cite[Theorem 3.3]{FKMT15}), this leads to the stability estimate
\begin{equation}\label{eq:stabatomic}
\int_{\R^d} \ip{\nu^1_{t,x}}{\big|\xi-v(x,t)\big|}\,dx \leq \int_{\R^d} \ip{\bar{\nu}^1_x}{\big|\xi-\bar{v}(x)\big|}\,dx
\end{equation}
for any entropy solution $v$. Thus, if $\bar{\nu}^1_x = \delta_{\bar{u}(x)}$ then also $\nu^1_{t,x} = \delta_{u(x,t)}$---in other words, \eqref{eq:entrcond1} provides stability with respect to entropy solutions $u(x,t)$, realized as atomic entropy measure-valued solutions. Note, however, that if $\bar{\nu}$ is non-atomic then the right-hand side of \eqref{eq:stabatomic} is $O(1)$. Hence, \eqref{eq:entrcond1} only imposes stability with respect to atomic statistical solutions. We propose instead the following:
\begin{center}
\emph{\textbf{Entropy condition:} The physically meaningful statistical solution must be stable not just with respect to\\ single constant functions, but to any finite convex combination of constant functions.}
\end{center}
Since constant functions do not lie in $L^1(\R^d)$, we need to introduce the following auxiliary lemma, which characterizes the set of transport plans, $\Pi(\mu,\rho)$, when $\rho$ is a convex combination of Dirac measures.
\begin{lemma}
Let $\mu,\rho\in\Prob(\F)$ such that $\rho$ is of the form $\rho = \sum_{i=1}^M \alpha_i \delta_{u_i}$ for coefficients $\alpha_i\geq0$, $\sum_i\alpha_i=1$ and functions $u_1,\dots,u_M\in\F$. Then a measure $\pi$ lies in $\Pi(\mu,\rho)$ if and only if there are $\mu_1,\dots,\mu_M\in\Prob(\F)$ such that
\begin{equation}\label{eq:transplanequiv}
\pi = \sum_{i=1}^M \alpha_i\mu_i\otimes\delta_{u_i}
\end{equation}
(and, in particular, $\sum_{i=1}^M \alpha_i\mu_i = \mu$).
\end{lemma}
\begin{proof}
Necessity is immediate. For sufficiency, let $\pi\in\Pi(\mu,\rho)$ and define $\mu_i(A) := \frac{\pi(A\times\{u_i\})}{\alpha_i}$. Without loss of generality, we may assume that $\alpha_i>0$ and that $u_1,\dots,u_M$ are distinct. Since $\pi(\F\times\{u_i\}) = \rho(\{u_i\}) = \alpha_i$ we have $\mu_i\in\Prob(\F)$ for each $i$. Moreover, $\pi(A\times\{u_i\}) = \alpha_i\mu_i(A) = \alpha_i(\mu_i\otimes\delta_{u_i})(A\times\{u_i\})$ for each $i$, so \eqref{eq:transplanequiv} follows.
\end{proof}
Based on this simple observation we conclude that whenever $\rho$ is $M$-atomic with weights $\alpha_i$, there is a one-to-one correspondence between transport plans $\pi\in\Pi(\mu,\rho)$ and elements of the set
\[
\Lambda(\alpha,\mu) := \Bigl\{(\mu_1,\dots,\mu_M)\ :\ \textstyle\sum_{i=1}^M \alpha_i\mu_i = \mu \Bigr\} \qquad \text{for } \alpha=(\alpha_1,\dots,\alpha_M),\ \alpha_i\geq0,\ \sum_{i=1}^M \alpha_i=1.
\]
This set is never empty since $(\mu,\dots,\mu) \in \Lambda(\alpha,\mu)$ for any $\alpha$. Note that the set $\Lambda(\alpha,\mu)$ depends on the target measure $\rho$ \emph{only} through the weights $\alpha_1,\dots,\alpha_M$.

\begin{definition}
A statistical solution $\mu_t$ is termed an \define{entropy statistical solution} if for every choice of coefficients $\alpha_i>0$ with $\sum_{i=1}^M \alpha_i=1$ and for every $(\bar{\mu}_1,\dots,\bar{\mu}_M)\in\Lambda(\alpha,\bar{\mu})$, there exists a map $t\mapsto(\mu_{1,t},\dots,\mu_{M,t})\in\Lambda(\alpha,\mu_t)$ such that $\mu_{i,0}=\bar{\mu}_i$ and
\begin{equation}\label{eq:entrcondss}
\begin{split}
\sum_{i=1}^M\alpha_i\left[\int_{\R_+}\int_{\F}\int_{\R^d} \big|u(x)-c_i\big|\partial_t\phi + q\big(u(x),c_i\big)\cdot \nabla_x\phi\,dxd\mu_{i,t}(u)dt +
\int_{\F}\int_{\R^d}\big|\bar{u}(x)-c_i\big|\phi\Bigr|_{t=0}\,dxd\bar{\mu}_i(\bar{u})\right]  \geq 0
\end{split}
\end{equation}
for all $0\leq\phi\in C_c^\infty(\R^d\times\R_+)$ and for all constants $c_1,\dots,c_M\in\R$. (Here, $q(u,c)$ is the Kruzkov entropy flux function.)
\end{definition}
\begin{lemma}
The canonical statistical solution is an entropy statistical solution.
\end{lemma}
\begin{proof}
%If $\mu_t=S_t\#\bar{\mu}$ then $\Lambda(\alpha,\mu_t) = \left\{(S_t\#\bar{\mu}_i)_{i=1}^M\ :\ (\mu_i)_{i=1}^M\in\Lambda(\alpha,\bar{\mu})\right\}$. 
Select $(\bar{\mu}_1,\dots,\bar{\mu}_M)\in\Lambda(\alpha,\bar{\mu})$ for an arbitrary weight $\alpha$ and define $\mu_{i,t} := S_t\#\bar{\mu}_i$. Then $(\mu_{1,t},\dots,\mu_{M,t})\in\Lambda(\alpha,\mu_t)$, and
\begin{align*}
&\sum_{i=1}^M\alpha_i\left[\int_{\R_+}\int_{\F}\int_{\R^d} \bigl|u(x)-c_i\bigr|\partial_t\phi + q\big(u(x),c_i\big)\cdot \nabla_x \phi\,dxd\mu_{i,t}(u)dt + \int_{\F}\int_{\R^d}\big|\bar{u}(x)-c_i\big|\phi\Bigr|_{t=0}\,dxd\bar{\mu}_i(\bar{u})\right] \\
=&\ \sum_{i=1}^M\alpha_i\int_\F\left[\int_{\R_+}\int_{\R^d} \big|S_t\bar{u}(x)-c_i\big|\partial_t\phi + q\big(S_t\bar{u}(x),c_i\big)\cdot \nabla_x \phi\,dxdt + \int_{\R^d}\big|\bar{u}(x)-c_i\big|\phi\Bigr|_{t=0}\,dx\right]d\bar{\mu}_i(\bar{u}) \\
\geq&\ 0,
\end{align*}
since the map $(x,t)\mapsto S_t\bar{u}(x)$ is an entropy solution of the deterministic problem.
\end{proof}
Note that 

\begin{lemma}\label{lem:convdiracstable}
Let $\mu_t$ be an arbitrary entropy statistical solution with initial data $\bar{\mu}\in\Prob(\F)$ satisfying $\supp\bar{\mu}\subset\EuScript{U}$. Fix $\alpha_1,\dots,\alpha_M>0$ with $\sum_{i=1}^{M}\alpha_i=1$. Let $v_1,\dots,v_M : \R_+\to\EuScript{U}$ be entropy solutions of \eqref{eq:cl} with initial data $\bar{v}_1,\dots,\bar{v}_M\in\EuScript{U}$, respectively, and define
\[
\bar{\rho}:=\sum_{i=1}^M\alpha_i\delta_{v_i}, \qquad \rho_t:=\sum_{i=1}^M\alpha_i\delta_{v_i(t)} \qquad \forall\ t\in\R_+.
\]
Then
\begin{equation}\label{eq:ssstab}
W_1(\rho_t,\mu_t) \leq W_1(\bar{\rho},\bar{\mu}) \qquad \forall\ t>0.
\end{equation}
\end{lemma}
\begin{proof}
Let $(\bar{\mu}_i)_{i=1}^M\in\Lambda(\alpha,\bar{\mu})$ be an optimal transport plan from $\bar{\mu}$ to $\bar{\rho}$. The entropy condition for $\mu_t$ gives the existence of a map $t\mapsto\big(\mu_{i,t}\big)_{i=1}^n$ such that
\begin{equation}\label{eq:kruzkov0}
\sum_{i=1}^M\alpha_i\left[\int_{\R_+}\int_{\F}\int_{\R^d} \big|u(x)-c_i\big|\partial_t\phi + q\big(u(x),c_i\big)\cdot \nabla_x \ \phi\,dxd\mu_{i,t}(u)dt + \int_{\F}\int_{\R^d}\big|\bar{u}(x)-c_i\big|\phi\Bigr|_{t=0}\,dxd\bar{\mu}_i(\bar{u})\right] \geq 0
\end{equation}
for any choice of $\phi\in C_c^\infty(\R^d \times\R_+)$ and $c_i\in\R$. Let $\phi = \phi(x,y,t,s)\in C_c^\infty((\R^d)^2\times\R_+^2)$. Set $c_i = v_i(y,s)$ for some point $(y,s)$ and integrate over $y\in\R$ and $s\in\R_+$:
\begin{equation}\label{eq:kruzkov1}
\begin{split}
\int_{\R_+}\int_{\R^d}\sum_{i=1}^M\alpha_i\Biggl[\int_{\R_+}\int_{\F}\int_{\R^d} \big|u(x)-v_i(y,s)\big|\partial_t\phi + q\big(u(x),v_i(y,s)\big)\cdot \nabla_x \phi\,dxd\mu_{i,t}(u)dt \\
+ \int_{\F}\int_{\R^d}\big|\bar{u}(x)-v_i(y,s)\big|\phi\Bigr|_{t=0}\,dxd\bar{\mu}_i(\bar{u})\Biggr]dyds \geq 0.
\end{split}
\end{equation}
(The expression in the brackets is measurable with respect to $(y,s)$ since \eqref{eq:kruzkov0} is continuous with respect to $c_i$.)

Next, since each $v_i$ is an entropy solution, we have for all $\xi\in\R$ and $0\leq\phi\in C_c^\infty(\R^d\times\R_+)$
\[
\int_{\R_+}\int_{\R^d} \big|\xi-v_i(y,s)\big|\partial_s\phi + q\big(\xi,v_i(y,s)\big)\cdot \nabla_y\phi\,dyds + \int_{\R^d} \big|\xi-\bar{v}_i(y)\big|\phi\Bigr|_{s=0}\,dy \geq 0.
\]
Set $\xi=u(x)$ for some $u\in\F$ and $x\in\R$. Integrate the above over $x\in\R$ and over $u\in\F$ with respect to $\mu_{i,t}$ for some $t\in\R_+$. Integrate over $t\in\R_+$, multiply by $\alpha_i$ and sum over $i=1,\dots,M$:
\begin{equation}\label{eq:kruzkov2}
\begin{split}
\sum_{i=1}^{M}\alpha_i\int_{\R_+}\int_\F\int_{\R^d} \Biggl[\int_{\R_+}\int_{\R^d} \big|u(x)-v_i(y,s)\big|\partial_s\phi + q\big(u(x),v_i(y,s)\big)\cdot \nabla_y\phi\,dyds \\
+ \int_{\R^d} \big|u(x)-\bar{v}_i(y)\big|\phi\Bigr|_{s=0}\,dy \Biggr]dxd\mu_{i,t}(u)dt \geq 0.
\end{split}
\end{equation}
Applying Fubini's theorem to this and equation \eqref{eq:kruzkov1} and adding the two, we obtain
\begin{gather*}
\sum_{i=1}^M\alpha_i \Biggl[\int_{\R_+}\int_{\R_+}\int_\F\int_{\R^d}\int_{\R^d} \big|u(x)-v_i(y,s)\big|(\partial_t+\partial_s)\phi + q\big(u(x),v_i(y,s)\big)\cdot (\nabla_x+\nabla_y)\phi\,dxdyd\mu_{i,t}(u)dtds \\
+ \int_{\R_+}\int_\F\int_{\R^d}\int_{\R^d} \big|u(x)-\bar{v}_i(y)\big|\phi\Bigr|_{s=0}\,dxdyd\mu_{i,t}(u)dt + \int_{\R_+}\int_{\F}\int_{\R^d}\int_{\R^d} \big|\bar{u}(x)-v_i(y,s)\big|\phi\Bigr|_{t=0}\,dxdyd\bar{\mu}_i(\bar{u})ds\Biggr] \geq 0.
\end{gather*}
Now set $\phi(x,y,t,s) := \psi\Bigl(\frac{x+y}{2}, \frac{t+s}{2}\Bigr)\omega_\eps(x-y)\omega_{\eps'}(t-s)$ for some nonnegative $\psi\in C_c^\infty(\R^d\times\R_+)$ and a mollifier $\omega_\eps$. Using the dominated convergence theorem on the integrals over $\F$, we find that as $\eps\to0$, the above converges to
\begin{gather*}
\sum_{i=1}^M\alpha_i \Biggl[\int_{\R_+}\int_{\R_+}\int_\F\int_{\R^d} \big|u(x)-v_i(x,s)\big|(\partial_t+\partial_s)\tilde{\phi} + 2q\big(u(x),v_i(x,s)\big)\cdot \nabla_x\tilde{\phi}\,dxd\mu_{i,t}(u)dtds \\
+ \int_{\R_+}\int_\F\int_{\R^d} \big|u(x)-\bar{v}_i(x)\big|\tilde{\phi}\Bigr|_{s=0}\,dxd\mu_{i,t}(u)dt + \int_{\R_+}\int_{\F}\int_{\R^d} \big|\bar{u}(x)-v_i(x,s)\big|\tilde{\phi}\Bigr|_{t=0}\,dxd\bar{\mu}_i(\bar{u})ds\Biggr] \geq 0,
\end{gather*}
where $\tilde{\phi}(x,t,s) := \psi\Bigl(x,\frac{t+s}{2}\Bigr)\omega_{\eps'}(t-s)$. Finally, letting $\eps'\to0$ we get
\begin{gather*}
\sum_{i=1}^M\alpha_i \Biggl[\int_{\R_+}\int_\F\int_{\R^d} \big|u(x)-v_i(x,t)\big|\partial_t\psi + q\big(u(x),v_i(x,t)\big)\cdot \nabla_x\psi\,dxd\mu_{i,t}(u)dt \\
+ \int_\F\int_{\R^d} \big|\bar{u}(x)-\bar{v}_i(x)\big|\psi\Bigr|_{t=0}\,dxd\bar{\mu}_i(\bar{u})\Biggr] \geq 0.
\end{gather*}
We now set $\psi(x,\tau) := \ind_{[0,t]}(\tau)$ for some $t\in\R_+$ to get
\[
\sum_{i=1}^M\alpha_i\Biggl[-\int_\F \big\|u-v_i(t)\big\|_{\F}\,d\mu_{i,t}(u) + \int_\F \big\|\bar{u}-\bar{v}_i\big\|_\F\,d\bar{\mu}_i(\bar{u})\Biggr] \geq 0.
\]
Using the fact that $(\bar{\mu}_i)$ is an optimal transport plan from $\bar{\mu}$ to $\bar{\rho}$, we end up with \eqref{eq:ssstab}.
\end{proof}

To complete our proof of well-posedness of statistical solutions we need the following well-known result, whose proof is included in the appendix for the sake of completeness.
\begin{lemma}\label{lem:diracdense}
Let $X$ be a Polish space equipped with its Borel $\sigma$-algebra. Then the convex hull of Dirac measures on $X$ is dense in $\Prob(X)$ with respect to the topology of weak convergence. In other words, for every $\mu\in\Prob(X)$, there is a sequence $\rho_n\in\Prob(X)$ of convex combinations of Dirac measures such that $\rho_n\wto\mu$ as $n\to\infty$.
\end{lemma}

\begin{theorem}
Let $\bar{\mu}\in\Prob(\F)$ with $\supp\bar{\mu}\subset\EuScript{U}:=L^1\cap L^\infty(\R^d)$. Then the entropy statistical solution with initial data $\bar{\mu}$ is unique and coincides with the canonical statistical solution. Any two entropy statistical solutions $\mu_t$, $\rho_t$ satisfy
\begin{equation}
W_1(\mu_t,\rho_t) \leq W_1(\bar{\mu},\bar{\rho}).
\end{equation}
\end{theorem}
\begin{proof}
Let $\mu_t$ be an entropy statistical solution with initial data $\bar{\mu}$. By Lemma \ref{lem:diracdense}, the convex hull of Dirac measures is dense in $\Prob(\F)$, so we can find a sequence $\bar{\mu}_n\in\Prob(\F)$ ($n\in\N$) of convex combinations of Dirac measures such that $\bar{\mu}_n \wto \bar{\mu}$ in $\Prob(\F)$ as $n\to\infty$. Let $\mu_{n,t} := S_t\#\bar{\mu}_n$ be the corresponding canonical statistical solutions, and note that also $\mu_{n,t} \wto S_t\#\bar{\mu}$ as $n\to\infty$. From Lemma \ref{lem:convdiracstable} we find that
\[
W_1(\mu_t,\mu_{n,t}) \leq W_1(\bar{\mu},\bar{\mu}_n) \to 0 \qquad \text{as } n\to\infty.
\]
Thus, $\mu_t = \wlim_{n\to\infty}\mu_{n,t} = S_t\#\bar{\mu}$, whence $\mu_t$ is the canonical statistical solution.
\end{proof}

\section{Discussion}
\label{sec:disc}
Given the lack of global in time existence results, and the recent non-uniqueness results of \cite{DLS1,CDL3}, the acceptance of entropy solutions as the standard solution paradigm for multi-dimensional systems of conservation laws is being increasingly questioned. Based on extensive numerical results, recent papers such as \cite{FKMT15} have advocated entropy measure-valued solutions (MVS), as defined by DiPerna \cite{DiP85}, as an appropriate solution paradigm for systems of conservation laws. However, entropy MVS are not necessarily unique, even for scalar conservation laws, if the MVS is non-atomic. Since numerical results of \cite{FKMT15} strongly hint at the possibility of non-atomic MVS even when the initial data is a atomic, it is natural to seek additional constraints on entropy MVS to enforce uniqueness.

Given this background, and the need for developing a solution concept that can accommodate uncertain initial data (and corresponding uncertain solutions) that arise frequently in the area of uncertainty quantification (UQ), we seek to adapt the notion of statistical solutions, originally developed in \cite{Foi72,Foi73} for the incompressible Navier--Stokes equations, to systems of conservation laws. Statistical solutions are time-parametrized probability measures on some (infinite-dimensional) function space. Infinite-dimensional Liouville or Hopf equations track the evolution of the time-parametrized measure. However, the extension of statistical solutions as defined in \cite{Foi72,Foi73,FMRT}, to systems of conservation laws, is highly non-trivial as the ``natural'' function spaces for the dynamics of conservation laws consists merely of integrable functions, and may lack the regularity required to define the Liouville or Hopf equations. Although one can work with probability measures on distributions in the specific case of the inviscid Burgers equation (as suggested in \cite{CD1,CD2,Bert}), it is very difficult to enforce uniqueness on such a large space of measures. Another disadvantage of probability measures on functions is that they do not readily provide any local (statistical) information at specific (collections of) points in the spatial domain.

We define statistical solutions for systems of conservation laws in a different manner. To this end, we prove a {novel} equivalence theorem between probability measures on $L^p$ spaces ($1 \leq p < \infty$) and a family (hierarchy) of Young measures, the so-called correlation measures, on finite-dimensional tensor product spatial domains. For all $k \in \N$, the $k$-th member of this hierarchy, the so-called $k$-point correlation marginal, is a Young measure that provides information on correlations of the underlying functions at $k$ distinct points in the spatial domain. In particular, the first correlation marginal is classical one-point Young measure. Thus, a probability measure on an $L^p$ space can be realized as an Young measure, augmented with multi-point correlations on the spatial domain. This representation enables us to \emph{localize} probability measures on function spaces and view them as a collection of all possible multi-point correlation marginals. We also show that moments of the correlation marginals uniquely determine the corresponding probability measure on the infinite-dimensional function space. We believe that this representation of probability measures will be of independent interest in stochastic analysis, particularly stochastic partial differential equations \cite{DZ}, in uncertainty quantification of evolutionary PDEs \cite{UQhb} and in Bayesian inversion and data assimilation for time-dependent PDEs \cite{ST1}. In particular, the use of statistical solutions will provide a framework for uncertainty quantification that does not depend on any particular parametrization of the solution in terms of random fields, as is customary in UQ \cite{UQhb}. 

In this paper, we use the equivalence between probability measures on $L^p$ and families of correlation measures to define statistical solutions of systems of conservation laws. In particular, we utilize the fact that moments of correlation measures uniquely determine the underlying probability measure, to evolve these moments in a manner consistent with the dynamics of the system \eqref{eq:cl}. Thus, a statistical solution has to satisfy an (infinite) family of nonlinear PDEs, but each of these PDEs is defined on a finite-dimensional (tensor-product) spatial domain. This should be contrasted with the infinite-dimensional Liouville or Hopf equations that the statistical solutions of \cite{Foi72,Foi73,FMRT} need to satisfy. Moreover, our notion of statistical solutions restricts the class of probability measures to those on $L^p$ spaces, rather than on distributions (as in \cite{CD2}) and makes it more amenable to analysis, particularly from the point of view of uniqueness. At the same time, our notion of statistical solutions augment the standard concept of measure-valued solutions, with additional information in the form of multi-point correlations, and paves the way for constraining the solutions sufficiently to guarantee uniqueness.

We investigate the well-posedness of the proposed concept of statistical solutions in the specific context of multi-dimensional scalar conservation laws in this paper. We show existence by proving that the push forward of the initial probability measure on $L^1 \cap L^{\infty}$ by the Kruzkhov entropy solution semi-group is a statistical solution, and we term this solution the \emph{canonical statistical solution}. We propose a novel admissibility criteria, based on stability with respect to a suitable stationary statistical solution, namely probability measures supported on finite collections of constant functions. These \emph{entropy statistical solutions} are a generalization of the standard Kruzkhov entropy solutions for scalar conservation laws. We show that the canonical statistical solution is the unique entropy statistical solution. Furthermore, we show that it is contractive with respect to the $1$-Wasserstein metric on probability measures on $L^1$. Thus, entropy statistical solutions for multi-dimensional scalar conservation laws are shown to be well-posed and are thus completely characterized. 

This article is the first in a series of papers investigating statistical solutions of multi-dimensional systems of conservation laws. We lay out the measure theoretic basis, define statistical solutions for systems and show well-posedness in the scalar case. Forthcoming papers in the series will deal with numerical approximation of entropy statistical solutions of scalar conservation laws \cite{FLyeSid1} and global existence of statistical solutions for a large class of multi-dimensional systems of conservation laws by showing convergence of a Monte Carlo based numerical approximation algorithm \cite{FKLyeSid1}. Admissibility criteria that single out physically relevant statistical solutions are the topic of current and future work.

\section*{Acknowledgments}
U.S.F.\ was supported in part by the grant \textit{Waves and Nonlinear Phenomena} (WaNP) from the Research Council of Norway. S.M. was supported in part by ERC STG. N 306279, SPARCCLE. The authors thanks Kjetil O. Lye and Franziska Weber (SAM, ETH) for their helpful comments.

%\appendix
%\section{Measure theory}
%\subsection{Kolmogorov's extension theorem}
%For Euclidean sets $U$ and $D$, we denote by $\U^D$ the space of all functions $u:D\to\U$. For $u\in\U^D$ and $x=(x_1,\dots,x_k)\in D^k$ we define the canonical projection $\pi^k_x:\U^D\to\U^k$ as $\pi^k_x(u) = (u(x_1),\dots,u(x_k))$. In particular, $\pi^1$ is just the identity mapping, $\pi^1(u) = u$. Let $\Salg:=\bigotimes_{x\in D}\Borel(\U)$ denote the product $\sigma$-algebra on $\U^D$ (see e.g.\ Klenke \cite[Chapter 14]{Kle}).
%
%\begin{theorem}[Kolmogorov's extension theorem]\label{thm:kolmext}
%For all finite collections of points $x_1,\dots,x_k\in D$ let $P_{x_1,\dots,x_k}\in\Prob(\U^k)$ satisfy
%\begin{itemize}
%\item $P_{x_{\sigma(1)},\dots,x_{\sigma(k)}}(F_1\times\dots\times F_k) = P_{x_1,\dots,x_k}\bigl(F_{\sigma^{-1}(1)}\times\dots\times F_{\sigma^{-1}(k)}\bigr)$ for all permutations $\sigma$ of $\{1,\dots,k\}$.
%\item $P_{x_1,\dots,x_{k+1}}(F_1\times\dots\times F_k\times\R) = P_{x_1,\dots,x_{k}}(F_1\times\dots\times F_k)$.
%\end{itemize}
%Then there exist a unique probability measure $P$ on $(\U^D,\,\Salg)$ such that for all $k\in\N$ and all distinct $x_1,\dots,x_k\in D$, we have $\pi^k_{x_1,\dots,x_k}\# P = P_{x_1,\dots,x_k}$.
%\end{theorem}
%\begin{proof}
%See \cite[Theorem 14.36]{Kle} or \cite[Theorem 2.1.5]{Oks}.
%\end{proof}

\appendix
\section{Appendix}
For completeness we provide the proof of Proposition \ref{prop:uniquecyl}. The proof relies on the following two lemmas.
\begin{lemma}\label{lem:cylring}
$\Cyl(X)$ is a ring.\footnote{A collection of sets $\EuScript{X}\subset 2^X$ is a \define{ring} if $\emptyset\in\EuScript{X}$ and if both $A\cup B$ and $A\setminus B$ lie in $\EuScript{X}$ whenever $A,B\in\EuScript{X}$.}
\end{lemma}
\begin{proof}
%Abbreviate $\lip{\phi}{u} = \phi(u)$. 
Clearly, $\emptyset\in\Cyl(X)$, and if $A_1, A_2\in\Cyl(X)$ are of the form
\[
A_i = \left\{u\in X\ :\ \bigl(\phi_1^i, \dots, \phi_{n_i}^i\bigr)(u) \in F_i\right\}, \qquad i=1,2
\]
then both
\[
A_1\cup A_2 = \left\{u\in X\ :\ \bigl(\phi_1^1, \dots, \phi_{n_1}^1, \phi_1^2, \dots, \phi_{n_2}^2\bigr)(u) \in \left(F_1\times \R^{n_2}\right) \cup \left(\R^{n_1}\times F_2\right)\right\}
\]
and
\[
A_1\setminus A_2 = \left\{u\in X\ :\ \bigl(\phi_1^1, \dots, \phi_{n_1}^1, \phi_1^2, \dots, \phi_{n_2}^2\bigr)(u) \in F_1\times \bigl(F_2\bigr)^c\right\}
\]
are cylinder sets.
\end{proof}

\begin{lemma}\label{lem:normidentity}
%If $X$ is a normed vector space that is either separable, or is the dual of a separable space, 
If $X$ is a separable normed vector space
then there exists a countable family $\{\phi_n\}_{n\in\N} \subset X^*$ such that
\begin{equation}\label{eq:normidentity}
\|u\|_X=\sup_{n\in\N}\,\lip{\phi_n}{u} \quad \text{ for every $u\in X$}.
\end{equation}
\end{lemma}
\begin{proof}
%If $X$ is separable,
Let $\{u_n\}_{n\in\N} \subset X$ be a countable dense subset of the unit sphere $\partial B_1(0) \subset X$. For each $n\in\N$, let $\phi_n\in X^*$ satisfy $\lip{\phi_n}{u_n} = 1$ and $\|\phi_n\|_{X^*} = 1$. If $u\in \partial B_1(0)$ is arbitrary and $\eps>0$, find an $u_n$ such that $\|u-u_n\|_X < \eps$. Then
\[
1 \geq \lip{\phi_n}{u} = \lip{\phi_n}{u_n} - \lip{\phi_n}{u_n-u} \geq 1-\eps,
\]
so $\|u\|_X=1$ can be approximated from below by $\lip{\phi_n}{u}$. Equation \eqref{eq:normidentity} follows.
%
%If $X = Y^*$ for a separable normed vector space $Y$, let $\{y_n\}_{n\in\N}\subset Y$ be a dense subset of the unit sphere in $Y$. Let $u\in X$. Since $X=Y^*$ we have $\|u\|_X = \sup_{y\in Y : \|y\|_Y=1} u(y)$, so if $y\in Y$ with $\|y\|_Y=1$ is such that $\|u\|_X - u(y) < \eps$ and $y_n$ is such that $\|y_n-y\|_Y < \eps$ then 
%\[
%0\leq \|u\|_X - u(y_n) < \eps + \eps\|u\|_X.
%\]
%Hence, the norm of $u$ can be approximated by $u(y_n)$. But now every $y_n$ can be identified with an element $\phi_n$ of the double dual $Y^{**} = X^*$ as $\phi_n(u) = u(y_n)$.
\end{proof}

\begin{proof}[Proof of Proposition \ref{prop:uniquecyl}]
Let $\{\phi_n\}_{n\in\N}$ be as in Lemma \ref{lem:normidentity}. For a $u_0\in X$ and $r>0$, the open ball of radius $r$ with centre $u_0$ can be written
\begin{align*}
B_r(u_0) &= \Bigl\{u\in X\ :\ \lip{\phi_n}{u-u_0} < r\ \forall\ n\in\N\Bigr\} \\
&= \bigcap_{n\in\N} \Bigl\{u\in X\ :\ \lip{\phi_n}{u} \in \bigl(-\infty,\ \lip{\phi_n}{u_0}+r\bigr)\Bigr\},
\end{align*}
which is a countable intersection of cylinder sets. It follows that $\sigma(\Cyl(X))$, the $\sigma$-algebra generated by $\Cyl(X)$, contains the $\sigma$-algebra generated by the open balls in $X$, which is precisely $\Borel(X)$. But every cylinder set is a Borel set; hence the two $\sigma$-algebras coincide, and \textit{(i)} follows.

By Lemma \ref{lem:cylring}, $\Cyl(X)$ is a ring which, by \textit{(i)}, generates $\Borel(X)$. Assertion \textit{(ii)} then follows from the fact that (signed) measures vanishing on a ring, vanish on the $\sigma$-algebra generated by the ring.
\end{proof}

\begin{proof}[Proof of Lemma \ref{lem:diracdense}]
Recall that the topology of weak convergence on $\Prob(X)$ for a Polish metric space $X$ is the coarsest topology for which the map $\mu\mapsto \int \phi \,d\mu$ is continuous for every $\phi\in C_b(X)$ \cite[Remark 13.14(ii)]{Kle}. Thus, the topology of weak convergence is generated by the open sets
\[
U_{\phi,\mu,\eps} := \left\{\rho\in\Prob(X)\ :\ \Bigl|\int\phi\, d\mu - \int\phi\,d\rho\Bigr| < \eps\right\}
\]
for $\mu\in\Prob(X)$, $\eps>0$ and $\phi\in C_b(X)$. It suffices to show that every nonempty open set $U_{\phi,\mu,\eps}$ contains a measure which is a convex combination of Dirac measures. Let $\bar{\phi}(x) = \sum_{i=1}^n a_i\ind_{A_i}(x)$ be a simple function such that $\sup_{x\in X}|\phi(x)-\bar{\phi}(x)| < \eps/2$. Fix $x_i\in A_i$ and define $\rho := \sum_{i=1}^n \mu(A_i)\delta_{x_i}$. Since $|\phi(x_i)-\phi(x)|<\eps$ for every $x\in A_i$, we find that
\begin{gather*}
\bigg|\int_X \phi\,d\rho - \int_X \phi\,d\mu\bigg| = \biggl|\sum_{i=1}^n\int_{A_i}\phi(x_i)- \phi(x)\,d\mu\biggr| \leq \sum_{i=1}^n\int_{A_i}|\phi(x_i)-\phi(x)|\,d\mu < \eps.
\end{gather*}
Hence, $\rho \in U_{\phi,\mu,\eps}$.
\end{proof}

%\bibliographystyle{plain}
%\bibliography{bibliography}

\begin{thebibliography}{10}

\bibitem{AGS05}
L.~Ambrosio, N.~Gigli, and G.~Savar{\'e}.
\newblock {\em {Gradient Flows: In Metric Spaces and in the Space of
  Probability Measures}}.
\newblock {Lectures in Mathematics ETH Z{\"u}rich}. Birkh{\"a}user Basel, 2005.

\bibitem{Bal89}
J.~Ball.
\newblock {A version of the fundamental theorem for {Young} measures}.
\newblock In M.~Rascle, D.~Serre, and M.~Slemrod, editors, {\em {PDEs and
  Continuum Models of Phase Transitions}}, volume 344 of {\em {Lecture Notes in
  Physics}}, pages 207--215. Springer Berlin / Heidelberg, 1989.
  
  \bibitem{BTW1}
C. Bardos, E. Titi and E. Wiedemann.
\newblock  The vanishing viscosity as a selection principle for the Euler equations: the case of 3D shear flow.
\newblock {\em C.R. Math. Acad. Sci. Paris}, 350 (15-16), 2012, 757-760.
  
\bibitem{SBG1}
S. Benzoni-Gavage and D. Serre.
\newblock Multidimensional hyperbolic partial differential equations. First order systems and applications. 
\newblock {\em Oxford university press}, 2007.

\bibitem{Bert}
J. Bertoin.
\newblock The inviscid Burgers equation with Brownian initial velocity.
\newblock {\em Comm. Math. Phys.,} 193 (2), 1998, 397-406.

\bibitem{BB1}
S. Bianchini and A. Bressan.
\newblock Vanishing viscosity solutions of nonlinear hyperbolic systems.
\newblock \textit{Ann. of Math.} (2) 161 (2005), no. 1, 223--342.

\bibitem{lncse:uq}
H. Bijl, D. Lucor, S. Mishra and Ch. Schwab. (editors).
\newblock {\em Uncertainty quantification in computational fluid dynamics.,}
\newblock Lecture notes in computational science and engineering 92, Springer, 2014. 

\bibitem{BRE1}
A. Bressan.
\newblock Hyperbolic systems of conservation laws: The one dimensional Cauchy problem.
\newblock {\em Oxford university press,} 2000.

\bibitem{CD1}
L. Carraro and J. Duchon.
\newblock Intrinsic statistical solutions of the Burgers equation and Levy processes.
\newblock {\em C.R. Math. Acad. Sci. Paris} 319 (8), 1994, 855-858.

\bibitem{CD2}
L. Carraro and J. Duchon.
\newblock Burgers equation with initial conditions with homogeneous and independent increments.
\newblock {\em Ann. Inst. H. Poincar\'e Anal. Non Lineare,} 15 (4). 1998, 431-458.

\bibitem{Chae1}
D. Chae,
\newblock The vanishing viscosity limit of statistical solutions of the Navier--Stokes equations. I. 2-D periodic case.
\newblock {\em J. Math. Anal. Appl.,} 155 (2), 1991, 437-459.

\bibitem{Chae2}
D. Chae,
\newblock The vanishing viscosity limit of statistical solutions of the Navier--Stokes equations. II. The general case.
\newblock {\em J. Math. Anal. Appl.,} 155 (2), 1991, 460-484.

\bibitem{CGLM08}
Pierre Comon, Gene Golub, Lek-Heng Lim, and Bernard Mourrain.
\newblock {Symmetric Tensors and Symmetric Tensor Rank}.
\newblock {\em SIAM Journal on Matrix Analysis and Applications},
  30(3):1254--1279, 2008.

\bibitem{DZ}
G.~{Da Prato} and J.~Zabczyk.
\newblock {\em {Stochastic equations in infinite dimensions}}.
\newblock Cambridge University Press, 1992.

\bibitem{DAF1}
C. Dafermos.
\newblock \textit{Hyperbolic conservation laws in continuum physics.}
\newblock Springer, Berlin, 2000.

\bibitem{DLS1}
C. De Lellis, L. Sz\'ekelyhidi Jr.
\newblock The Euler equations as a differential inclusion.
\newblock \textit{Ann. of Math.} (2) 170 (2009), no. 3, 1417--1436.

\bibitem{CDL3}
E. Chiodaroli, C. De Lellis, O. Kreml.
\newblock Global ill-posedness of the isentropic system of gas dynamics.
\newblock {\em Comm. Pure Appl. Math.,} 68 (7), 2015, 1157-1190.

\bibitem{DST12}
S. Demoulini, D. M. A. Stuart and A. E. Tzavaras.
\newblock Weak-strong uniqueness of dissipative measure-valued solutions for polyconvex elastodynamics.
\newblock \textit{Archive for Rational Mechanics and Analysis} 205(3), 927--961, 2012.

\bibitem{DU77}
J.~Diestel and J.~J. Uhl.
\newblock {\em {Vector Measures}}.
\newblock American Mathematical Society, 1977.

\bibitem{DiP85}
R.~J. DiPerna.
\newblock {Measure-valued solutions to conservation laws}.
\newblock {\em Archive for Rational Mechanics and Analysis}, 88:223--270, 1985.

\bibitem{DM87}
R. J. DiPerna and A. Majda.
\newblock Oscillations and concentrations in weak solutions of the incompressible fluid equations.
\newblock \textit{Comm. Math. Phys.} 108 (4), 1987, 667--689.

\bibitem{Edw}
R.~E. Edwards.
\newblock {\em {Functional Analysis. Theory and Applications}}.
\newblock Holt, Rinehart and Winston, Inc., 1965.

\bibitem{FKMT15}
U.~S. Fjordholm, R.~K{\"a}ppeli, S.~Mishra, and E.~Tadmor.
\newblock {Construction of approximate entropy measure-valued solutions for
  hyperbolic systems of conservation laws}.
\newblock \textit{J. FoCM}, to appear, 2016, available from doi:10.1007/s10208-015-9299-z

\bibitem{FMTacta}
U. S. Fjordholm, S. Mishra and E. Tadmor.
\newblock On the computation of measure-valued solutions.
\newblock {\em Acta Numerica}, 2016, to appear. 

\bibitem{FLyeSid1}
U. S. Fjordholm, K. O. Lye and S. Mishra.
\newblock Statistical solutions of hyperbolic conservation laws II: Numerical approximation in the scalar case.
\newblock {\em In preparation,} 2016.

\bibitem{FKLyeSid1}
U. S. Fjordholm, R. K\"appeli, K. O. Lye and S. Mishra.
\newblock Statistical solutions of hyperbolic conservation laws III: Numerical approximation for multi-dimensional systems.
\newblock {\em In preparation,} 2016.

\bibitem{Foi72}
C. Foia\c s.
\newblock Statistical study of Navier--Stokes equations I.
\newblock \textit{Rend. Sem. Mat. Univ. Padova} 48, 219--348, 1972.

\bibitem{Foi73}
C. Foia\c s.
\newblock Statistical study of Navier--Stokes equations II.
\newblock \textit{Rend. Sem. Mat. Univ. Padova} 49, 9-123, 1973.

\bibitem{FMRT}
C. Foia\c s, O. Manley, R. Rosa, R. Temam.
\newblock \textit{Navier--Stokes Equations and Turbulence}. Cambridge University Press, 2001.

\bibitem{Frisch}
U. Frisch.
\newblock \textit{Turbulence}, Cambridge University Press, 1995.

\bibitem{UQhb}
R. Ghanem, D. Higdon and H. Owhadi (eds).
\newblock {\em Handbook of uncertainty quantification}, Springer, 2016.

\bibitem{GL1}
J. Glimm.
\newblock Solutions in the large for nonlinear hyperbolic systems of equations.
\newblock {\em Comm. Pure Appl. Math.}  18 (4), 1965, 697-715.

\bibitem{GGW1}
P. Gwiazda, A. Swierczewska-Gwiazda and E. Wiedemann.
\newblock Weak-Strong
uniqueness for measure-valued solutions of some compressible 
fluid models.
\newblock {\em Nonlinearity} 28, 2015, 3873-3890.

\bibitem{HR}
H.~Holden and N.~H. Risebro.
\newblock {\em {Front Tracking for Hyperbolic Conservation Laws}}.
\newblock Springer-Verlag Berlin Heidelberg, 2011.

\bibitem{IL1}
R. Illner and J. Wick.
\newblock On statistical and measure-valued solutions of differential equations.
\newblock {\em J. Math. Anal. Appl.,} 157 (2), 1991, 351-365.

\bibitem{Iss18}
L. Isserlis
\newblock On a formula for the product-moment coefficient of any order of a normal frequency distribution in any number of variables.
\newblock \textit{Biometrika} 12: 134--139, 1918.

\bibitem{Kle}
A.~Klenke.
\newblock {\em {Probability Theory. A Comprehensive Course}}.
\newblock Springer London, 2nd edition, 2014.

\bibitem{Kruz1}
S.~N. Kruzkov.
\newblock {First order quasilinear equations in several independent variables}.
\newblock {\em Math USSR SB}, 10(2):217--243, 1970.

\bibitem{GL2}
H. Lim, Y. Yu, J. Glimm, X. L. Li and D. H. Sharp.
\newblock Chaos, transport and mesh convergence for fluid mixing.
\newblock {\em Act. Math. Appl. Sin.,} 24 (3), 2008, 355--368.

\bibitem{Pan1}
E. Yu. Panov.
\newblock On the statistical solutions of the Cauchy problem for a first-order quasilinear equation (Russian).
\newblock {\em Mat. Model.,} 14 (3), 2002, 17-26.

\bibitem{Sch89}
S. Schochet.
\newblock Examples of measure-valued solutions.
\newblock \textit{Communications in Partial Differential Equations} 14(5), 545--575, 1989.

\bibitem{ST1}
A. M. Stuart.
\newblock Inverse problems: a Bayesian perspective.
\newblock {\em Acta Numerica,} 19, 2010, 451-559.

\bibitem{Vil}
C. Villani.
\newblock {\em {Topics in Optimal Transportation}}.
\newblock {Graduate Studies in Mathematics, Vol. 58}. American Mathematical Society, 2003.

\bibitem{Øks}
B. \O ksendal.
\newblock \textit{Stochastic Differential Equations. An Introduction with Applications},
\newblock 6th edition, Springer Berlin Heidelberg, 2003.
\end{thebibliography}

\end{document}